\documentclass[12pt,a4paper]{article}
\usepackage{mathrsfs}
\usepackage{amssymb}
\usepackage{amsmath}
\usepackage{amsfonts}
\usepackage{amstext}
\usepackage{amsthm}
\usepackage{graphicx}
\usepackage[top=1in, bottom=1in,left=0.75in,right=0.75in]{geometry}

\def \qed {\hfill \vrule height6pt width 6pt depth 0pt}

\begin{document}

\title{Eigenvalues, bifurcation and one-sign solutions for the periodic $p$-Laplacian
\thanks{Research supported by the NSFC (No.11061030, No.10971087).}}
\author{{\small  Guowei Dai$^{a}$\thanks{Corresponding author. \newline
\text{\quad\,\,\, E-mail address}: daiguowei@nwnu.edu.cn (G. Dai), mary@nwnu.edu.cn (R. Ma), Haiyan.Wang@asu.edu (H. Wang).%
},\, \, Haiyan Wang$^{a,b}$
} \\
{\small $^{a}$Department of Mathematics, Northwest Normal University, Lanzhou, 730070, PR China}\\
{\small $^{b}$School of Mathematical and Natural Sciences}\\
{\small Arizona State University, Phoenix, AZ 85069, USA}\\
}
\date{}
\maketitle

\begin{abstract}
In this paper, we establish a unilateral global bifurcation result for a class of quasilinear periodic boundary problems with a sign-changing weight.
By the Ljusternik-Schnirelmann theory, we first study the spectrum of the periodic $p$-Laplacian with the sign-changing weight.
In particular, we show that there exist two simple, isolated, principal eigenvalues $\lambda_0^+$ and $\lambda_0^-$.
Furthermore, under some natural hypotheses on perturbation function,
we show that $\left(\lambda_0^\nu,0\right)$ is a bifurcation
point of the above problems and there are two distinct unbounded sub-continua
$\mathscr{C}_\nu^{+}$ and $\mathscr{C}_\nu^{-}$,
consisting of the continuum $\mathscr{C}_\nu$ emanating from $\left(\lambda_0^\nu, 0\right)$,
where $\nu\in\{+,-\}$. As an application of the above result, we study the existence of one-sign solutions
for a class of quasilinear periodic boundary problems with the sign-changing weight.
Moreover, the uniqueness of one-sign solutions and the dependence of solutions on the parameter $\lambda$ are
also studied.
\end{abstract}\textbf{\ }

\numberwithin{equation}{section}

\numberwithin{equation}{section}

\section{Introduction}

\quad\, In the past few decades, periodic boundary value problems have attracted the attention of many specialists in differential equations because of
their interesting applications. For example, the application in looking for spatially
periodic solutions of the well-known Camassa--Holm equation, see [\ref{ACHM}, \ref{CH}, \ref{CHH}, \ref{C}, \ref{CM}].
The Camassa--Holm equation is a recently discovered model for
the propagation of shallow water waves of moderate amplitude [\ref{CL}, \ref{J}] and some authors have already indicated recently that the
equation might be relevant to the modeling of tsunamis [\ref{CJ}, \ref{La}]. As the recent year examples,
we mention the papers of Atici and Guseinov [\ref{AG}], Jiang et al. [\ref{JCOA}], Li [\ref{L}],
O'Regan and Wang [\ref{OW}], Torres [\ref{T}], Zhang and Wang [\ref{ZW}], Graef et al. [\ref{GKW}] and references therein.
Their main tool is the fixed-point theorem of cone expansion/compression type. Ma et al. [\ref{MXH}, \ref{MXH1}]
studied the existence of positive solutions for the second-order periodic boundary value problems by making use of the bifurcation techniques.

Recently, Dai and Ma [\ref{DM}] established unilateral global bifurcation theory
for one-dimensional $p$-Laplacian problems with 0-Dirichlet boundary condition.
Moreover, Dai and Ma [\ref{DM}], Dai [\ref{D}] also studied the existence of nodal solutions
for the one-dimensional $p$-Laplacian problems based on the unilateral global bifurcation theory.
For the abstract unilateral global bifurcation theory, we refer the reader to [\ref{DM}, \ref{D1}, \ref{D2}, \ref{L1}]
and the references therein.

The main purpose of this paper is to establish a result similar to that of [\ref{DM}]
about the continuum of one-sign solutions for the following periodic $p$-Laplacian problem
\begin{equation}\label{pb}
\left\{
\begin{array}{l}
-\left(\varphi_p\left(u'\right)\right)'+q(x)\varphi_p(u)=\lambda m(x)\varphi_p(u)+g(x,u,\lambda),\,\,0<x<T,\\
u(0)=u(T),\,\, u'(0)=u'(T),
\end{array}
\right.
\end{equation}
where $1<p<+\infty$, $\varphi_p(s)=\vert s\vert^{p-2}s$, 
$q \in C([0,T]; [0,+\infty))$ with $q\not\equiv 0$, $m:[0,T]\rightarrow \mathbb{R}$ is a sign-changing
weight and $g:[0,T]\times \mathbb{R}^2\rightarrow\mathbb{R}$ is continuous satisfying $g(x,s,0)\equiv 0$.
We also assume that the perturbation function $g$ satisfies the following hypothesis
\begin{equation}\label{c1}
\lim_{s\rightarrow0}\frac{g(x,s,\lambda)}{\vert s\vert^{p-1}}=0
\end{equation}
uniformly on $[0,T]$ and $\lambda$ on bounded sets.

In order to study the unilateral global bifurcation phenomena of problem (\ref{pb}), we must consider the
following eigenvalue problem
\begin{equation}\label{pe}
\left\{
\begin{array}{l}
-\left(\varphi_p\left(u'\right)\right)'+q(x)\varphi_p(u)=\lambda m(x)\varphi_p(u),\,\,0<x<T,\\
u(0)=u(T),\,\, u'(0)=u'(T).
\end{array}
\right.
\end{equation}
For the case of $p=2$, Constantin [\ref{C1}] has proved that problem (\ref{pe}) possesses
two infinite sequences of eigenvalues
\begin{equation}
\cdots<\lambda_2^-\leq\lambda_1^-<\lambda_0^-<\lambda_0^+<\lambda_1^+\leq \lambda_2^+<\cdots\nonumber
\end{equation}
such that $\lambda_0^+$ and $\lambda_0^-$ are simple eigenvalues with positive eigenfunctions.
However, the methods used in [\ref{C1}] cannot be used to deal with problem (\ref{pe}) because
$p$-Laplace operator is neither self-adjoint linear nor symmetric.
For $m(x)\equiv 1$ or $m(x)>0$ on $[0,T]$, by using the variational method, Binding and Rynne [\ref{BR1}, \ref{BR2}] have shown that
problem (\ref{pe}) has a sequences of eigenvalues
\begin{equation}
-\infty<\lambda_0<\lambda_0^\nu\leq\cdots.\nonumber
\end{equation}
Moreover, $\lambda_0$ is a simple, principal eigenvalue. Note that it is unnecessary for $\lambda_0>0$ because
$q$ is not necessarily positive in [\ref{BR1}, \ref{BR2}]. In this paper,
we also use the Ljusternik-Schnirelmann theory to study problem (\ref{pe}).
More precisely, we shall show that problem (\ref{pe})
possesses two infinite sequences of eigenvalues
\begin{equation}
\cdots\leq\lambda_2^-\leq\lambda_1^-<\lambda_0^-<0<\lambda_0^+<\lambda_1^+\leq \lambda_2^+\leq\cdots\nonumber
\end{equation}
and $\lambda_0^+$ and $\lambda_0^-$ are simple, isolated, principal eigenvalues and continuous with respect to $p$.
This method is also used by Cuesta [\ref{Cu}] to study eigenvalue problems for the $p$-Laplacian with 0-Dirichlet boundary condition and sign-changing weight.

Following the above spectrum results, we shall show that
$\left(\lambda_0^\nu,0\right)$ is a bifurcation point of one-sign solutions to problem (\ref{pb}) and there are two distinct unbounded sub-continua
$\mathscr{C}_\nu^+$ and $\mathscr{C}_\nu^-$, consisting of the continuum $\mathscr{C}_\nu$ bifurcating from $\left(\lambda_0^\nu, 0\right)$,
where $\nu\in\{+,-\}$.

On the basis of the unilateral global bifurcation result, we
investigate the existence of one-sign solutions for the following periodic $p$-Laplacian problem
\begin{equation}\label{pn}
\left\{
\begin{array}{l}
-\left(\varphi_p\left(u'\right)\right)'+q(x)\varphi_p(u)=\lambda m(x)f(u),\,\,0<x<T,\\
u(0)=u(T),\,\, u'(0)=u'(T),
\end{array}
\right.
\end{equation}
where $f\in C(\mathbb{R})$, $\lambda$ is a parameter.
Here, we shall establish some results of existence, multiplicity and nonexistence of one-sign solutions for
problem (\ref{pn}) according to the asymptotic behavior of $f$ at 0 and $\infty$ and the fact of whether $f$ possesses zeros in $\mathbb{R} \backslash\{0\}$.
Our results extend and improve the corresponding ones of [\ref{GKW}].
To the best of our knowledge, most results of this paper are new even in the case of $p=2$.
We now give a brief description of the contents of the paper.

In Section 2, with the aid of the Ljusternik-Schnirelmann theory and operator theory, we
study the variational eigenvalues of problem (\ref{pe}).
Moreover, as a byproduct, we also establish several important
properties of a quasilinear operator which itself possesses
an independent importance. The results of this section partially extend the corresponding ones of [\ref{BR1}, \ref{BR2}].

In Section 3, we prove some properties of the principle eigenvalues $\lambda_0^+$ and $\lambda_0^-$. More precisely,
we shall show that $\lambda_0^+$ and $\lambda_0^-$ are simple, isolated, principal eigenvalues (their corresponding
eigenfunctions are positive or negative) and continuous with respect to $p$. It is well-known that the continuity of
$\lambda_0^+$ and $\lambda_0^-$ with respect to $p$ is crucial in the studying of the global bifurcation phenomena for
$p$-Laplacian. We use the method established by Del Pino et al. [\ref{DPEM}, \ref{DPM}] to prove this result but with some
extra effort since the boundary condition is different from [\ref{DPM}].
 To the best of our knowledge, this result is new even in the case of $m\geq 0$.

In Section 4, we establish the unilateral global bifurcation theory for problem (\ref{pb}).
In the global bifurcation theory of differential equations,
it is well-known that a change of the index of the trivial solution implies
the existence of a branch of nontrivial solutions, bifurcating from the set of trivial solutions
which is either unbounded or returns to the set of trivial solution.
Hence, the index formula of an isolated zero is very important in
the study of the bifurcation phenomena for differential equations.
Firstly, we establish an index formula for $p=2$ by the linear compact operator theory.
Then by use of the index formula and the deformation along $p$, we prove an index formula
involving the problem (\ref{pe}) which guarantees $\left(\lambda_0^\nu,0\right)$ is a bifurcation point of nontrivial solutions
to problem (\ref{pb}). Furthermore, by an argument similar to that of [\ref{DM}], we can get unilateral
global bifurcation results for problem (\ref{pb}).

In Section 5, we study the existence of one-sign solutions for problem (\ref{pn}) with signum condition according
to the asymptotic behavior of $f$ at 0 and $\infty$. The results of this section extend and improve the corresponding ones of
[\ref{GKW}, Theorem 2.1] and [\ref{MXH1}, Theorem 1.1 and 1.2] even in the case of $p=2$.

In Section 6, we show a result involving the uniqueness and dependence of solutions on the parameter.
This result extends and improves the corresponding ones to [\ref{GKW}, Theorem 2.2] even in the case of $p=2$.
To prove this result, we introduce a new method which is different from that of [\ref{GKW}, \ref{LL}, \ref{LL1}].

Finally, Section 7 is devoted to study the existence of one-sign solutions for problem (\ref{pn}) without signum condition.
To do this, following some ideas from [\ref{R0}], we establish a unilateral global bifurcation theorem from infinity
for problem (\ref{pb}). This theorem, as an independent result, is of interest too. Our results of this section extend and improve the corresponding results of [\ref{MXH}].

\section{Variational eigenvalues}

\quad\, In this section, we shall establish the eigenvalue theory for problem (\ref{pe}) via the Ljusternik-Schnirelmann theory. Let
\begin{equation}
W_{T}^{1,p}(0,T):=\left\{u\in W^{1,p}(0,T)|u(0)=u(T)\right\}\nonumber
\end{equation}
with the norm
\begin{equation}
\Vert u\Vert=\left(\int_0^{T}\left(\left\vert u'\right\vert^p+q\vert u\vert^p\right)\,dx\right)^{\frac{1}{p}}.\nonumber
\end{equation}

It is not difficult to verify that $W_{T}^{1,p}(0,T)$ is a real Banach space. For simplicity, we write $u_n\rightharpoonup u$ and
$u_n\rightarrow u$ to indicate the weak
convergence and strong convergence of sequence $\left\{u_n\right\}$ in $W_{T}^{1,p}(0,T)$, respectively.\\

First, we recall the definition of weak solution.\\ \\
\textbf{Definition 2.1.} $u\in W_T^{1,p}(0,T)$ is called a weak solution of problem (\ref{pe}) if
\begin{equation}\label{ws}
\int_0^{T}\left(\left\vert u'\right\vert^{p-2}u'\phi'+q\vert u\vert^{p-2}u\phi\right)\,dx=\lambda\int_0^{T} m\vert u\vert^{p-2}u\phi\,dx\nonumber
\end{equation}
for any $\phi\in W_T^{1,p}(0,T)$.\\ \\
\indent For the regularity of weak solution, we have the following result.
\\ \\
\noindent\textbf{Lemma 2.1.} \emph{Any weak solution $u\in W_T^{1,p}(0,T)$ of problem (\ref{pe}) is also a classical solution of problem (\ref{pe}).}
\\ \\
\indent In order to prove Lemma 2.1, we need the following technical result.
\\ \\
\noindent\textbf{Proposition 2.1.} \emph{Let $f:\mathbb{R}\rightarrow \mathbb{R}$ be a
function. For a given $x_0\in \mathbb{R}$, if $f$
is continuous in some neighborhood $U$ of $x_0$, differential in $U\setminus \{x_0\}$
and $\underset{x\rightarrow x_0}\lim f'(x)$ exists,
then $f$ is differential at $x_0$ and $f'\left(x_0\right)=\underset{x\rightarrow x_0}\lim f'(x)$.}\\ \\
\textbf{Proof.} The conclusion is a direct corollary of Lagrange mean theorem, we omit the proof here.\qed
\\ \\
\textbf{Proof of Lemma 2.1.} According to Definition 2.1, we have
\begin{equation}
-\left(\left\vert u'\right\vert^{p-2}u'\right)'+q\vert u\vert^{p-2}u
=\lambda m\vert u\vert^{p-2}u\,\,\text{in}\,\,(0,T)\nonumber
\end{equation}
in the sense of distribution, i.e.,
\begin{equation}
-\left(\left\vert u'\right\vert^{p-2}u'\right)'+q\vert u\vert^{p-2}u
=\lambda m\vert u\vert^{p-2}u\,\,\text{in}\,\,(0,T)\setminus I\nonumber
\end{equation}
for some $I\subset (0,T)$ which satisfies $\text{meas}\{I\}=0$.
Clearly, the embedding of $W_T^{1,p}(0,T)\hookrightarrow
C^{\alpha}[0,T]$ with some $\alpha\in(0,T)$ is compact since $W_T^{1,p}(0,T)\hookrightarrow W^{1,p}(0,T)$
is continuous and $W^{1,p}(0,T)\hookrightarrow C^{\alpha}[0,T]$ is compact (see [\ref{E}]). Consequently, we obtain
\begin{equation}
-\left(\left\vert u'\right\vert^{p-2}u'\right)'\in C([0,T]\setminus I).\nonumber
\end{equation}
Set $v:=\varphi_p\left(u'\right)$. The above relation implies that
$\lim_{x\rightarrow x_0} v'$ exists for any $x_0\in I$. Hence,
Proposition 2.1 implies that $v\in C^1[0,T]$. By appropriate choosing of the test function $\phi$,
we can show that $u$ satisfies the first equation of problem (\ref{pe}).
Furthermore, using Definition 2.1 and integrating by parts, we can see that $u$ satisfies
the periodic boundary condition $u(0)=u(T)$, $u'(0)=u'(T)$.\qed\\

Define the functional on $W_T^{1,p}(0,T)$
\begin{equation}
\Phi(u)=\int_0^{T}\frac{1}{p}\left(\left\vert u'\right\vert^p+q\vert u\vert^p\right)\,dx.\nonumber
\end{equation}
It is obvious that the functional $\Phi$ is continuously G\^{a}teaux differentiable.
Denote $L:= \Phi': W_T^{1,p}(0,T) \rightarrow \left(W_T^{1,p}(0,T)\right)^*$; then
\begin{equation}\label{oL}
\langle L(u),v\rangle=\int_0^{T} \left(\left\vert u'\right\vert^{p-2}u'v'+q\vert u\vert^{p-2}uv\right)\,dx, \,\, \forall u, v\in W_T^{1,p}(0,T),
\end{equation}
where $\left(W_T^{1,p}(0,T)\right)^*$ denotes the dual space of $W_T^{1,p}(0,T)$; $\langle\cdot,\cdot\rangle$ is the duality pairing between
$W_T^{1,p}(0,T)$ and $\left(W_T^{1,p}(0,T)\right)^*$.
\\

We have the following properties about the operator $L$.
\\ \\
\textbf{Proposition 2.2.} \emph{(i) $L: W_T^{1,p}(0,T)\rightarrow \left(W_T^{1,p}(0,T)\right)^*$ is a continuous and strictly monotone operator;}

\emph{(ii) $L$ is a map of type $\left(S_+\right)$, i.e., if $u_n\rightharpoonup u$
in $W_T^{1,p}(0,T)$ and
\begin{equation}
\underset{n\rightarrow+\infty}{\overline{\lim}}
\left\langle L\left(u_n\right)-L(u),u_n-u\right\rangle\leq0,\nonumber
\end{equation}
then $u_n \rightarrow u$ in $W_T^{1,p}(0,T)$;}

\emph{(iii) $L:W_T^{1,p}(0,T)\rightarrow \left(W_T^{1,p}(0,T)\right)^*$ is a homeomorphism.}\\
\\
\textbf{Proof.}
(i) It is not difficult to verify that $L$ is continuous.
For any $u$, $v\in W_T^{1,p}(0,T)$ with $u\neq v$ in $W_T^{1,p}(0,T)$.
By the Cauchy's inequality, we have
\begin{equation}\label{ci}
uv\leq \vert u\vert \vert v\vert\leq\frac{\vert u\vert^2+\vert v\vert^2}{2}.
\end{equation}
Noting (\ref{ci}), we can easily obtain that
\begin{equation}\label{ci1}
\int_{0}^{T}\left\vert u'\right\vert^{p}\,{d}x-\int_{0}^{T}\left\vert u'\right\vert^{p-2} u' v'\,{d}x
\geq\int_{0}^{T}\frac{\left\vert u'\right\vert^{p-2}}{2}\left(\left\vert u'\right\vert^{2}-\left\vert  v'\right\vert^{2}\right)\,{d}x,
\end{equation}
\begin{equation}\label{ci11}
\int_{0}^{T}\vert u\vert^{p}\,{d}x-\int_{0}^{T}\vert u\vert^{p-2} u v\,{d}x
\geq\int_{0}^{T}\frac{\vert u\vert^{p-2}}{2}\left(\vert u\vert^{2}-\vert  v\vert^{2}\right)\,{d}x,
\end{equation}
\begin{equation}\label{ci2}
\int_{0}^{T}\left\vert v'\right\vert^{p}\,{d}x-\int_{0}^{T}\left\vert v'\right\vert^{p-2} u'v'\,{d}x
\geq\int_{0}^{T}\frac{\left\vert v'\right\vert^{p-2}}{2}\left(\left\vert v'\right\vert^{2}-\left\vert u'\right\vert^{2}\right)\,{d}x,
\end{equation}
and
\begin{equation}\label{ci22}
\int_{0}^{T}\vert v\vert^{p}\,{d}x-\int_{0}^{T}\vert v\vert^{p-2} uv\,{d}x
\geq\int_{0}^{T}\frac{\vert v\vert^{p-2}}{2}\left(\vert v\vert^{2}-\vert u\vert^{2}\right)\,{d}x.
\end{equation}
By virtue of (\ref{oL}), (\ref{ci1}),  (\ref{ci11}), (\ref{ci2}) and (\ref{ci22}), we obtain that
\begin{eqnarray}\label{oLm0}
\langle L(u)-L(v),u-v\rangle&=&\langle L(u),u\rangle-\langle L(u),v\rangle-\langle L(v),u\rangle+\langle L(v),v\rangle\nonumber\\
&=&\left(\int_{0}^{T}\left(\left\vert u'\right\vert^{p}+q\vert u\vert^p\right)\,{d}x
-\int_{0}^{T}\left(\left\vert u'\right\vert^{p-2}u'v'+q\vert u\vert^{p-2} uv\right)\,{d}x\right)\nonumber\\
& &-\left(\int_{0}^{T}\left(\left\vert v'\right\vert^{p-2}v'u'+q\vert v\vert^{p-2}vu\right)\,{d}x-\int_{0}^{T}\left(\left\vert v'\right\vert^{p}+q\vert v\vert^p\right)\,{d}x\right)\nonumber\\
&\geq&\int_{0}^{T}\frac{\left\vert u'\right\vert^{p-2}}{2}\left(\left\vert u'\right\vert^{2}-\left\vert v'\right\vert^{2}\right)\,{d}x-
\int_{0}^{T}\frac{\left\vert v'\right\vert^{p-2}}{2}\left(\left\vert u'\right\vert^{2}-\left\vert v'\right\vert^{2}\right)\,{d}x\nonumber\\
& &+\int_{0}^{T}\frac{\vert u\vert^{p-2}}{2}q\left(\vert u\vert^{2}-\vert v\vert^{2}\right)\,{d}x-
\int_{0}^{T}\frac{\vert v\vert^{p-2}}{2}q\left(\vert u\vert^{2}-\vert v\vert^{2}\right)\,{d}x\nonumber\\
&\geq&\int_{0}^{T}\frac{1}{2}
\left(\left\vert u'\right\vert^{p-2}-\left\vert v'\right\vert^{p-2}\right)\left(\left\vert u'\right\vert^{2}-\left\vert v'\right\vert^{2}\right)\,{d}x\nonumber\\
& &+\int_{0}^{T}\frac{1}{2}
q\left(\vert u\vert^{p-2}-\vert v\vert^{p-2}\right)\left(\vert u\vert^{2}-\vert v\vert^{2}\right)\,{d}x\geq 0,
\end{eqnarray}
i.e., $L$ is monotone. In fact, $L$ is strictly monotone. Indeed, if $\langle L(u)-L(v),u-v\rangle=0$, then we have
\begin{equation}
\left\vert u'\right\vert=\left\vert v'\right\vert\,\,\text{and}\,\, \vert u\vert=\vert v\vert.\nonumber
\end{equation}
Thus, we obtain
\begin{eqnarray}\label{oLm}
\langle L(u)-L(v),u-v\rangle&=&\langle L(u),u-v\rangle-\langle L(v),u-v\rangle\nonumber\\
&=&\int_{0}^{T} \left\vert u'\right\vert^{p-2}\left(u'-v'\right)^2\,dx+\int_{0}^{T} q\vert u\vert^{p-2}\left(u-v\right)^2\,dx\nonumber\\
&=&0.
\end{eqnarray}
If $1<p<2$, (\ref{oLm}) implies that $u'=v'$ and $u=v$, which is a contradiction. If $p\geq2$, (\ref{oLm}) implies that
$u'=v'$ and $u=v$ which contradicts $u\neq v$ in $W_T^{1,p}(0,T)$ or $\left\vert u'\right\vert\equiv 0\equiv\vert u\vert$.
If the later case occurs, we get $v=u\equiv0$, which is a contradiction.
Therefore, $\langle L(u)-L(v),u-v\rangle>0$. It follows that $L$ is a strictly monotone operator on $W_T^{1,p}(0,T)$.

(ii) From (i), if $u_n\rightharpoonup u$ and $\underset{n\rightarrow+\infty}{\overline{\lim}}
\left\langle L\left(u_n\right)-L(u),u_n-u\right\rangle\leq0$, then
\begin{equation}
\underset{n\rightarrow+\infty}{\lim}
\left\langle L\left(u_n\right)-L(u),u_n-u\right\rangle=0.\nonumber
\end{equation}
In view of (\ref{oLm0}), $u_n'$ ($u_n$) converges in measure to $u'$ ($u$) in $(0,T)$, so
we get a subsequence (which we still denote by $u_n$) satisfying $u_n'(x)\rightarrow u'(x)$ and $u_n(x)\rightarrow u(x)$, a.e.
$x\in (0,T)$. By Fatou's Lemma we get
\begin{equation}\label{FL}
\underset{n\rightarrow+\infty}{\underline{\lim}}\int_0^{T}
\frac{1}{p}\left(\left\vert u_n'\right\vert ^{p}+q\left\vert u_n\right\vert^p\right)\,{d}x\geq\int_0^{T} \frac{1}{p}
\left(\left\vert u'\right\vert ^{p}+\vert u\vert^p\right)\,{d}x.
\end{equation}
From $u_n\rightharpoonup u$ we have
$\underset{n\rightarrow+\infty}{\lim}
\left\langle L\left(u_n\right),u_n-u\right\rangle=\underset{n\rightarrow+\infty}{\lim}
\left\langle L\left(u_n\right)-L(u),u_n-u\right\rangle=0$.
On the other hand, by Young's inequality, we have
\begin{eqnarray}\label{YI}
\left\langle L\left(u_n\right),u_n-u\right\rangle&=&\int_0^{T} \left(\left\vert u_n'\right\vert^{p}
+q\left\vert u_n\right\vert^p\right)\,dx-\int_0^{T} \left(\left\vert u_n'\right\vert^{p-2} u_n'u'+q\left\vert u_n\right\vert^{p-2}u_n u\right)\,dx\nonumber\\
&\geq& \int_0^{T} \left\vert u_n'\right\vert^{p}\,dx-\int_0^{T} \left\vert u_n'\right\vert^{p-1}\left\vert u'\right\vert\,dx\nonumber\\
& &+\int_0^{T} q\left\vert u_n\right\vert^{p}\,dx-\int_0^{T} q\left\vert u_n\right\vert^{p-1}\vert u\vert\,dx\nonumber\\
&\geq&\int_0^{T} \frac{1}{p}\left(\left\vert u_n'\right\vert^{p}+q\left\vert u_n\right\vert^p\right)\,dx-\int_0^{T}\frac{1}{p}\left(\left\vert u'\right\vert^{p}+q\vert u\vert^p\right)\,dx.
\end{eqnarray}
According to (\ref{FL}) and (\ref{YI}) we obtain
\begin{equation}
\underset{n\rightarrow+\infty}{\lim}\int_0^{T} \frac{1}{p}
\left(\left\vert u_n'\right\vert^{p}+q\vert u_n\vert^p\right)\,{d}x=\int_0^{T} \frac{1}{p}\left(\left\vert u'\right\vert^{p}+q\vert u\vert^p\right)\,{d}x.\nonumber
\end{equation}
By a similar method to prove [\ref{FZ}, Theorem 3.1], we have
\begin{equation}
\underset{n\rightarrow+\infty}{\lim}\int_0^{T}\left(\left\vert u_n'-u'\right\vert^{p}+q\left\vert u_n-u\right\vert^p\right)\,{d}x=0.\nonumber
\end{equation}
Therefore, $u_n\rightarrow u$, i.e., $L$ is of type $\left(S_+\right)$.

(iii) It is clear that $L$ is an injection since $L$ is a strictly monotone operator on $W_T^{1,p}(0,T)$. Since
\begin{equation}
\lim_{\Vert u\Vert\rightarrow+\infty}\frac{\langle L(u),u\rangle}{\Vert u\Vert}
=\lim_{\Vert u\Vert\rightarrow+\infty}\frac{
\int_0^{T}\left(\left\vert u'\right\vert^{p}+q\vert u\vert^p\right)\,dx}{\Vert u\Vert}=+\infty,\nonumber
\end{equation}
$L$ is coercive, thus $L$ is a surjection in view of Minty-Browder Theorem (see [\ref{Z},
Theorem 26A]). Hence $L$ has an inverse map $L^{-1}:\left(W_T^{1,p}(0,T)\right)^*\rightarrow W_T^{1,p}(0,T)$. Therefore, the continuity
of $L^{-1}$ is sufficient to ensure $L$ to be a homeomorphism.

If $f_n$, $f\in \left(W_T^{1,p}(0,T)\right)^*$, $f_n\rightarrow f$, let $u_n =L^{-1}\left(f_n\right)$, $u = L^{-1}(f)$,
then $L\left(u_n\right) = f_n$, $L(u) = f$.
The coercive property of $L$ implies that $\left\{u_n\right\}$ is bounded in $W_T^{1,p}(0,T)$.
We can assume that $u_{n_k}\rightharpoonup u_0$ in $W_T^{1,p}(0,T)$.
By $f_{n_k} \rightarrow f$ in $\left(W_T^{1,p}(0,T)\right)^*$, we have
\begin{equation}
\lim_{k\rightarrow+\infty}\left\langle L\left(u_{n_k}\right)-L\left(u_0\right),u_{n_k}-u_0\right\rangle
=\lim_{n\rightarrow+\infty}\left\langle f_{n_k}-f,u_{n_k}-u_0\right\rangle=0.\nonumber
\end{equation}
Since $L$ is of type $\left(S_+\right)$, $u_{n_k}\rightarrow u_0$.
Furthermore, the continuity of $L$ implies that $L\left(u_0\right)=L(u)$.
By injectivity of $L$, we have $u_0 = u$. So $u_{n_k}\rightarrow u$.
We claim that $u_n\rightarrow u$ in $W_T^{1,p}(0,T)$. Otherwise, there would exist a
subsequence $\left\{u_{m_j}\right\}$ of $\left\{u_n\right\}$ in $W_T^{1,p}(0,T)$ and an $\varepsilon_0 > 0$, such that for any $j\in \mathbb{N}$,
we have $\left\Vert u_{m_j}-u\right\Vert\geq \varepsilon_0$. But reasoning as above, $\left\{u_{m_j}\right\}$ would contain a
further subsequence $u_{m_{j_l}}\rightarrow u$ in $W_T^{1,p}(0,T)$ as $l \rightarrow+\infty$, which is a contradiction to
$\left\Vert u_{m_{j_l}}-u\right\Vert\geq \varepsilon_0$. Therefore, $L^{-1}$ is continuous.
\qed\\

Define the functional $\Psi:W_T^{1,p}(0,T)\rightarrow \mathbb{R}$ by
\begin{equation}
\int_0^{T}\frac{1}{p}m\vert u\vert^p\,dx.\nonumber
\end{equation}
\indent The following theorem is the main result of this section.
\\ \\
\textbf{Theorem 2.1.} \emph{The eigenvalue problem (\ref{pe}) has a sequence of eigenvalues
\begin{equation}
\cdots\leq\lambda_2^-\leq\lambda_1^-\leq\lambda_0^-<0<\lambda_0^+\leq\lambda_1^+\leq \lambda_2^+\leq\cdots.\nonumber
\end{equation}
Moreover,
\begin{equation}
\lambda_{0}^+=\inf\left\{\int_0^{T}\left(\left\vert u'\right\vert^p+q\vert u\vert^p\right)\,dx\Big| u\in W_T^{1,p}(0,T), \int_0^{T}m\vert u\vert^p\,dx=1\right\}\nonumber
\end{equation}
and}
\begin{equation}
\lambda_{0}^-=\max\left\{-\int_0^{T}\left(\left\vert u'\right\vert^p+q\vert u\vert^p\right)\,dx\Big| u\in W_T^{1,p}(0,T), \int_0^{T}-m\vert u\vert^p\,dx=1\right\}.\nonumber
\end{equation}
\noindent\textbf{Proof.} Set $\mathcal{M}=\{u\in W_T^{1,p}(0,T)\big|p\Psi(u)=1\}$ and
\begin{equation}
\Gamma_k=\{K\subset \mathcal{M}\big|K\text{\,\,is symmetric, compact and\,\,}\gamma(K)\geq k\}, \nonumber
\end{equation}
where $\gamma(K)$ is the genus of $K$. Then the weak form (also classical form by Lemma 2.1) of problem (\ref{pe}) on $\mathcal{M}$ can be equivalently written as
\begin{equation}\label{fz0}
\Phi'(u)=\lambda\Psi'(u),\,\, u\in \mathcal{M}.
\end{equation}
It is clear that ($\lambda,u$) solves (\ref{fz0}) if and only if $u$ is a critical point of $\Phi$ with respect to $\mathcal{M}$.
It is easy to verify that $\mathcal{M}$
is a closed symmetric $C^1$-submanifold of $W_T^{1,p}(0,T)$ with $0\not\in \mathcal{M}$, and ${\Phi}\in C^1({\mathcal{M},\mathbb{R}})$ is even.
It is obvious that ${\Phi}$ is bounded from below.

We claim that ${\Phi}$ satisfies the Palais-Smale condition at any level set $c$.

Suppose that $\{u_n\}\subset \mathcal{M}$, $\left\vert
\Phi\left(u_n\right)\right\vert\leq c$ and $ \Phi'\left(u_n\right)\rightarrow 0.$ Then for any constant $\theta>p$, we get
\begin{eqnarray}
c+\left\Vert u_{n}\right\Vert&\geq&{\Phi}\left(u_n\right)-\frac{1}{\theta}\Phi'\left(u_n\right)u_n \nonumber\\
&\geq&\left(\frac{1}{p}-\frac{1}{\theta}\right)\int_0^{T}\left(\left\vert u_n'\right\vert^{p}+q\vert u_n\vert^p\right)\,{d}x\nonumber\\
&=&\left(\frac{1}{p}-\frac{1}{\theta}\right)\left\Vert u_n\right\Vert^p.\nonumber
\end{eqnarray}
Hence, $\left\{\left\Vert u_n\right\Vert\right\}$ is bounded. Up to a subsequence,
we may assume that $u_n \rightharpoonup u$ in $\mathcal{M}$, so
$\left\langle\Phi'\left(u_n\right)-\Phi'(u),u_n-u\right\rangle\rightarrow 0$. By Proposition 2.2 (ii), we have $u_n\rightarrow u$.
Obviously, Proposition 2.2 (iii) implies that $0$ is not the eigenvalue of problem (\ref{pe}).
Now, applying Corollary 4.1 of [\ref{S}],
we obtain that problem (\ref{pe}) possesses a sequence of positive eigenvalues
\begin{equation}
0<\lambda_0^+\leq\lambda_1^+\leq \lambda_2^+\leq\cdots.\nonumber
\end{equation}
Moreover, we have that
\begin{equation}
\lambda_k^+=\inf_{K\in \Gamma_{k+1}}\sup_{u\in K}p\Phi(u).\nonumber
\end{equation}
In particular, if $k=0$, taking $K=\{u,-u| u\in\mathcal{M}\}$, we can get that
\begin{equation}
\lambda_{0}^+=\inf_{u\in \mathcal{M}}p\Phi(u)=\inf\left\{\int_0^{T}\left(\left\vert u'\right\vert^p
+q\vert u\vert^p\right)\,dx\Big| u\in W_T^{1,p}(0,T), \int_0^{T}m\vert u\vert^p\,dx=1\right\}.\nonumber
\end{equation}

In the case of $\lambda<0$, we restate eigenvalue problem (\ref{pe}) as the following
\begin{equation}\label{pe1}
\left\{
\begin{array}{l}
-\left(\varphi_p\left(u'\right)\right)'+q(x)\varphi_p(u)=\widehat{\lambda} \widehat{m}(x)\varphi_p(u),\,\,0<x<T,\\
u(0)=u(T),\,\, u'(0)=u'(T),
\end{array}
\right.
\end{equation}
where $\widehat{\lambda}=-\lambda$, $\widehat{m}(x)=-m(x)$. Using the above result, we have
that (\ref{pe1}) possesses a sequence of positive eigenvalues
\begin{equation}
0<\widehat{\lambda}_0^+\leq\widehat{\lambda}_1^+\leq \widehat{\lambda}_2^+\leq\cdots.\nonumber
\end{equation}
Set
\begin{equation}
\lambda_k^-:=-\widehat{\lambda}_k^+\nonumber
\end{equation}
for all $k\geq 0$. Thus, problem (\ref{pe}) also possesses a sequence of negative eigenvalues
\begin{equation}
\cdots\leq\lambda_2^-\leq\lambda_1^-\leq\lambda_0^-<0.\nonumber
\end{equation}
Similar to $\lambda_0^+$, we also get that
\begin{equation}
\lambda_{0}^-=\max\left\{-\int_0^{T}\left(\left\vert u'\right\vert^p+q\vert u\vert^p\right)\,dx\Big| u\in W_T^{1,p}(0,T),
\int_0^{T}-m\vert u\vert^p\,dx=1\right\}.\nonumber
\end{equation}
This completes the proof.\qed\\
\\
\textbf{Remark 2.1.} For $\nu\in\{+,-\}$ and each $k\geq 0$, Lemma 2.1 implies that $\lambda_k^\nu$ is a (classical) eigenvalue of problem (\ref{pe}).\\
\\
\textbf{Remark 2.2.} Note that if $m\geq0$ but $m\not\equiv 0$, we can only get the positive eigenvalues.

\section{Properties of positive minimal and negative maximal eigenvalues}

\quad\, In this section, we are going to study the properties of $\lambda_0^+$ and $\lambda_0^-$.
These properties, such as simplicity, isolation and the continuity with respect to $p$, are
important in the study of the global bifurcation phenomena for
$p$-Laplace problems, see [\ref{DPM}, \ref{LS}, \ref{P}].
\\

Similar to the results of the positive weight [\ref{BR2}, Theorem 3.1], we have the following theorem.\\ \\
\textbf{Theorem 3.1.} \emph{The eigenvalues $\lambda_0^+$ and $\lambda_0^-$ have the following properties.}\\

\emph{1. If $\lambda_0^-<\lambda<\lambda_0^+$ then problem (\ref{pe}) has no nontrivial solution.}

\emph{2. The eigenfunctions associated to $\lambda_0^+$ or $\lambda_0^-$ are either positive or negative on $[0,T]$.}

\emph{3. $\lambda_0^+$ and $\lambda_0^-$ are simple in the sense that the eigenfunctions associated to them are merely a constant multiple of each other.}

\emph{4. Any eigenfunction $u$ associated to $\lambda\neq \lambda_0^+$ and $\lambda\neq\lambda_0^-$ changes sign.}
\\ \\
\textbf{Proof.} We only consider the case of $\lambda\geq0$ since the proof of $\lambda<0$ can be given similarly.
Using a proof similar to that of [\ref{BR2}, Theorem 3.1] with obvious changes, we can obtain the properties of 1, 2 and 3.
However, the method which is used to prove Theorem 3.1 (c) of [\ref{BR2}] cannot be used directly here to prove 4 because $m$ is a sign-changing function.

Suppose on the contrary that $\lambda>\lambda_0^+$ and there exists an eigenfunction $u\geq 0$, i.e., $(\lambda,u)$ satisfies problem (\ref{pe}).
Similar to the proof of [\ref{BR2}, Theorem 3.1], we can show that $u>0$ on $[0,T]$.
Multiplying the first equation of problem (\ref{pe}) by $u$, we obtain after integration by parts
\begin{equation}
\int_0^{T}\left(\left\vert u'\right\vert^p+q\vert u\vert^p\right)\,dx=\lambda\int_0^{T}m\vert u\vert^p\,dx,\nonumber
\end{equation}
which implies that
\begin{equation}
\int_0^{T}m\vert u\vert^p\,dx>0.\nonumber
\end{equation}
So by scaling we may suppose that
\begin{equation}
\int_0^{T}m\vert u\vert^p\,dx>\frac{\lambda_0^+}{\lambda}.\nonumber
\end{equation}
Let $u_0^+$ be the eigenfunction corresponding to $\lambda_0^+$ satisfying $\int_0^{T}m\left\vert u_0^+\right\vert\,dx=1$.
Lemma 3.4 of [\ref{BR2}] yields
\begin{equation}
0\leq I\left(u,u_0^+\right)=\lambda_0^+-\lambda\int_0^{T}m u^p\,dx<0.\nonumber
\end{equation}
This is a contradiction.
\qed\\

Theorem 3.1 has shown that $\lambda_0^+$ is left-isolated and $\lambda_0^-$ is right-isolated. Furthermore,
we can show that $\lambda_0^+$ and $\lambda_0^-$ are isolated as the following.\\
\\
\textbf{Proposition 3.1.} \emph{$\lambda_0^+$ and $\lambda_0^-$ are isolated, that is, there exist $\delta^+>\lambda_0^+$ and $\delta^-<\lambda_0^-$
such that in the interval $\left(\delta^-,\delta^+\right)$ there is no other eigenvalues of problem (\ref{pe}).}
\\ \\
\textbf{Proof.} We only prove the isolated property of $\lambda_0^+$ since the case $\lambda_0^-$ is completely analogous.
Assume by contradiction that there exists a sequence of eigenvalues
$\lambda_n\in\left(\lambda_0^+, \delta^+\right)$ which converges to $\lambda_0^+$. Let $u_n$ be the corresponding eigenfunctions.
Theorem 3.1 implies that $u_n$ changes sign.
Integration by parts helps to yield
\begin{equation}
\int_0^{T}\left(\left\vert u_n'\right\vert^p+q\left\vert u_n\right\vert^p\right)\,dx=\lambda_n\int_0^{T}m\left\vert u_n\right\vert^p\,dx.\nonumber
\end{equation}
Define
\begin{equation}
v_n:=\frac{u_n}{\left(\int_0^{T} m(x)\left\vert u_n\right\vert^p\,dx\right)^{\frac{1}{p}}}.\nonumber
\end{equation}
Obviously, $v_n$ is bounded in $W_T^{1,p}(0,T)$ so there exists a subsequence, denoted again by $v_n$, and $v\in W_T^{1,p}(0,T)$ such that
$v_n\rightharpoonup v$ in $W_T^{1,p}(0,T)$ and $v_n\rightarrow v$ in $C[0,T]$.
Since functional $\Phi$ is sequentially weakly lower semi-continuous, we have
\begin{equation}
\int_0^{T}\left(\left\vert v'\right\vert^p+q\vert v\vert^p\right)\,dx \leq\liminf_{n\rightarrow+\infty}\int_0^{T}\left(\left\vert v_n'\right\vert^p+q\left\vert v_n\right\vert^p\right)\,dx=\liminf_{n\rightarrow+\infty}\lambda_n=\lambda_0^+.\nonumber
\end{equation}
On the other hand, $\int_0^{T} m\left\vert v_n\right\vert^p\,dx=1$ and $v_n\rightarrow v$ in $C[0,T]$
(hence, $\left\vert v_n\right\vert\rightarrow \vert v\vert$ uniformly on $[0,T]$) imply that
$\int_0^{T} m\vert v\vert^p\,dx=1$. It follows that
\begin{equation}
\int_0^{T}\left(\left\vert v'\right\vert^p+q\vert v\vert^p\right)\,dx \leq\lambda_0^+\int_0^{T} m\vert v\vert^p\,dx.\nonumber
\end{equation}
The above inequality and the variational characterization of $\lambda_0^+$ imply that
\begin{equation}
\int_0^{T} \left(\left\vert v'\right\vert^p+q\vert v\vert^p\right)\,dx=\lambda_0^+.\nonumber
\end{equation}
Then Theorem 3.1 follows that $v$ is positive or negative. Without loss of generality,
we may assume that $v>0$ on $[0,T]$. For any $\varepsilon\in\left(0,\min_{[0,T]} v\right)$, there exists $N_0>0$
such that $v_n\geq v-\varepsilon>0$ for any $n>N_0$ and all $x\in [0,T]$.
Thus $u_n\geq 0$ for $n$ large enough. This contradicts $u_n$ changing sign.\qed
\\

It is easy to see from Theorem 2.1 that the values of $\lambda_0^+$ and $\lambda_0^-$ are dependent on $p$. Hence,
we can rewrite $\lambda_0^+$ and $\lambda_0^-$ as $\lambda_0^+(p)$ and $\lambda_0^-(p)$ to indicate this dependence.
In fact, we can describe this relation more precisely as the following proposition does, and this proposition is crucial
to prove our main results in this paper.
\\ \\
\textbf{Proposition 3.2.} \emph{The eigenvalues $\lambda_0^+$ and $\lambda_0^-$ as the functions of $p$ are continuous.}
\\ \\
\indent In order to prove this proposition, we need the following results. \\
\\
\textbf{Lemma 3.1.} \emph{Let $C_T^\infty(\mathbb{R})$ be the space of indefinitely differentiable $T$-periodic functions
from $\mathbb{R}$ into $\mathbb{R}$. Then $C_T^\infty(\mathbb{R})$ is dense in $W_{T}^{1,p}(0,T)$.}
\\ \\
\textbf{Proof.} We define another norm on $W_{T}^{1,p}(0,T)$ by
\begin{equation}
\Vert u\Vert_*=\left(\int_0^{T}\left(\left\vert u'\right\vert^p+\vert u\vert^p\right)\,dx\right)^{\frac{1}{p}}.\nonumber
\end{equation}
It is easy to verify that $\Vert \cdot\Vert_*$ is equivalent to $\Vert \cdot\Vert$. From now on, we use $W_{T,0}^{1,p}(0,T)$ to denote
the space $W_{T}^{1,p}(0,T)$ with the norm $\Vert \cdot\Vert_*$.

In [\ref{MW}], Mawhin and Willem gave another definition of weak derivative which called $T$-weak derivative by Fan and Fan [\ref{FF}].
Let $\dot{u}$ denote the $T$-weak derivative of $u\in L^1(0,T)$. Define
\begin{equation}
\widetilde{W}_T^{1,p}(0,T)=\left\{u\in L^p(0,T)\big|\dot{u}\in L^p(0,T)\right\}\nonumber
\end{equation}
with the norm $\Vert u\Vert^*=\left(\int_0^T\left(\left\vert \dot{u}\right\vert^p+\vert u\vert^p\right)\,dx\right)^{1/p}$.
We also define $H_T^{1,p}(0,T)$ as the closure of $C_T^\infty(\mathbb{R})$ in $W^{1,p}(0,T)$. Lemma 2.3 of [\ref{FF}]
has shown that $H_T^{1,p}(0,T)=\widetilde{W}_T^{1,p}(0,T)$. It follows that $\Vert u\Vert^*=\Vert u\Vert_*$ for any
$u\in\widetilde{W}_T^{1,p}(0,T)$. Applying a similar method to prove [\ref{FF}, Lemma 2.11] with obvious changes,
we can show that $\widetilde{W}_T^{1,p}(0,T)=W_{T,0}^{1,p}(0,T)$. Hence, $C_T^\infty(\mathbb{R})$ is dense in
$W_{T,0}^{1,p}(0,T)$. Furthermore, $C_T^\infty(\mathbb{R})$ is dense in $W_{T}^{1,p}(0,T)$.\qed\\
\\
\textbf{Lemma 3.2.} \emph{We have}
\begin{equation}
\lambda_0^+(p)=\sup\left\{\lambda>0\Big|\lambda\int_0^{T}m\vert u\vert^p\,dx\leq\int_0^{T}\left(\left\vert u'\right\vert^p+q\vert u\vert^p\right)\,dx
\,\,\text{for all\,\,}u\in C_{T}^\infty(\mathbb{R})\right\}.\nonumber
\end{equation}
\noindent\textbf{Proof.} From the variational characterization of $\lambda_0^+(p)$ it follows that
\begin{equation}\label{cp21}
\lambda_0^+(p)=\sup\left\{\lambda>0\Big|\lambda\int_0^{T}m\vert u\vert^p\,dx\leq\int_0^{T}\left(\left\vert u'\right\vert^p+q\vert u\vert^p\right)\,dx
\,\,\text{for all\,\,}u\in W_{T}^{1,p}(0,T)\right\}.
\end{equation}
Let
\begin{equation}\label{cp22}
\lambda_0^\infty(p):=\sup\left\{\lambda>0\Big|\lambda\int_0^{T}m\vert u\vert^p\,dx\leq\int_0^{T}\left(\left\vert u'\right\vert^p+q\vert u\vert^p\right)\,dx
\,\,\text{for all\,\,}u\in C_{T}^\infty(\mathbb{R})\right\}.
\end{equation}

We claim that $\lambda_0^+(p)=\lambda_0^\infty(p)$.

Clearly, $\lambda_0^\infty(p)\leq \lambda_0^+(p)$. Next, we show that $\lambda_0^\infty(p)\geq \lambda_0^+(p)$.
Choosing $u\in W_T^{1,p}(0,T)$ but $u\neq cu_0^+$ for any $c\in \mathbb{R}$, Lemma 3.1 implies that there exists a sequence of $u_n\in C_{T}^\infty(\mathbb{R})$
such that $u_n\rightarrow u$ as $n\rightarrow +\infty$, where $u_0^+$ is the eigenfunction corresponding to $\lambda_0^+(p)$ with $\left\Vert u_0^+\right\Vert=1$.
There exists $n_0\in \mathbb{N}$ such that
\begin{equation}\label{cp24}
\lambda_0^+(p)\int_0^{T}m\vert u_{n_0}\vert^p\,dx\leq\int_0^{T}\left(\left\vert u_{n_0}'\right\vert^p+q\vert u_{n_0}\vert^p\right)\,dx.
\end{equation}
Otherwise, for all $n\in \mathbb{N}$, we have
\begin{equation}
\lambda_0^+(p)\int_0^{T}m\vert u_{n}\vert^p\,dx>\int_0^{T}\left(\left\vert u_{n}'\right\vert^p+q\vert u_{n}\vert^p\right)\,dx.\nonumber
\end{equation}
The fact that $X\hookrightarrow C[0,T]$ is compact and $u_n\rightarrow u$ in $X$ imply that
\begin{equation}\label{cp25}
\lambda_0^+(p)\int_0^{T}m\vert u\vert^p\,dx\geq\int_0^{T}\left(\left\vert u'\right\vert^p+q\vert u\vert^p\right)\,dx.
\end{equation}
While, (\ref{cp21}) implies that
\begin{equation}\label{cp23}
\lambda_0^+(p)\int_0^{T}m\vert u\vert^p\,dx\leq\int_0^{T}\left(\left\vert u'\right\vert^p+q\vert u\vert^p\right)\,dx.
\end{equation}
(\ref{cp25}) and (\ref{cp23}) imply that $u=cu_0^+$. This is a contradiction. Thus, (\ref{cp22}) and (\ref{cp24}) implies that $\lambda_0^+(p)\leq \lambda_0^\infty(p)$.
Therefore, we have
\begin{equation}
\lambda_0^+(p)=\sup\left\{\lambda>0\Big|\lambda\int_0^{T}m\vert u\vert^p\,dx\leq\int_0^{T}\left(\left\vert u'\right\vert^p+q\vert u\vert^p\right)\,dx
\,\,\text{for all\,\,}u\in C_{T}^\infty(\mathbb{R})\right\}.\nonumber
\end{equation}
\noindent\textbf{Proof of Proposition 3.2.} We only show that $\lambda_0^+:(1,+\infty)\rightarrow \mathbb{R}$ is continuous since the proof that $\lambda_0^-$ is
similar. In the following proof, we shall simply write  $\lambda_0^+$ as $\lambda_0$.

Lemma 3.2 has shown that
\begin{equation}\label{cp1}
\lambda_0(p)=\sup\left\{\lambda>0\Big|\lambda\int_0^{T}m\vert u\vert^p\,dx\leq\int_0^{T}\left(\left\vert u'\right\vert^p+q\vert u\vert^p\right)\,dx
\,\,\text{for all\,\,}u\in C_{T}^\infty(\mathbb{R})\right\}.
\end{equation}

Let $\{p_j\}_{j=1}^\infty$ be a sequence in $(1, +\infty)$ converging to $p > 1$. We shall show
that
\begin{equation}\label{cp2}
\lim_{j\rightarrow+\infty}\lambda_0\left(p_j\right)=\lambda_0(p).
\end{equation}
To do this, let $u\in C_{T}^\infty(\mathbb{R})$. Then, from (\ref{cp1}), we have that
\begin{equation}
\lambda_0\left(p_j\right)\int_0^{T}m\vert u\vert^{p_j}\,dx\leq\int_0^{T}\left(\left\vert u'\right\vert^{p_j}+q\vert u\vert^{p_j}\right)\,dx.\nonumber
\end{equation}
Applying the Dominated Convergence Theorem we find
\begin{equation}\label{cp3}
\limsup_{j\rightarrow+\infty}\lambda_0\left(p_j\right)\int_0^{T}m\vert u\vert^p\,dx\leq\int_0^{T}\left(\left\vert u'\right\vert^p+q\vert u\vert^p\right)\,dx.
\end{equation}
Relation (\ref{cp3}), the fact that $u$ is arbitrary and (\ref{cp1}) yield
\begin{equation}
\limsup_{j\rightarrow+\infty}\lambda_0\left(p_j\right)\leq\lambda_0(p).\nonumber
\end{equation}

Thus, to prove (\ref{cp2}) it suffices to show that
\begin{equation}\label{cp4}
\liminf_{j\rightarrow+\infty}\lambda_0\left(p_j\right)\geq\lambda_0(p).
\end{equation}
Let $\{p_k\}_{k=1}^\infty$ be a subsequence of $\{p_j\}_{j=1}^\infty$ such that
$\underset{k\rightarrow+\infty}\lim\lambda_0\left(p_k\right)=\underset{j\rightarrow+\infty}\liminf\lambda_0\left(p_j\right)$.

Let us fix $\varepsilon_0>0$ so that $p-\varepsilon_0>1$ and for each $0<\varepsilon<\varepsilon_0$,
$p-\varepsilon<p_k<p+\varepsilon$ if $k$ is large enough. For $k\in \mathbb{N}$,  choose $u_k\in W_{T}^{1,p_k}(0,T)$ such that
\begin{equation}\label{cp5}
\int_0^T\left(\left\vert u_k'\right\vert^{p_k}+q\left\vert u_k\right\vert^{p_k}\right)\,dx=1
\end{equation}
and
\begin{equation}\label{cp6}
\int_0^T\left(\left\vert u_k'\right\vert^{p_k}+q\left\vert u_k\right\vert^{p_k}\right)\,dx=\lambda_0\left(p_k\right)\int_0^T m\left\vert u_k\right\vert^{p_k}\,dx.
\end{equation}
For $0<\varepsilon<\varepsilon_0$ and $k$ large enough, (\ref{cp5}) and H\"{o}lder's inequality imply that
\begin{equation}\label{cp7}
\int_0^T\left\vert u_k'\right\vert^{p-\varepsilon}\,dx
\leq \left(1-\int_0^Tq\left\vert u_k\right\vert^{p_k}\,dx\right)^{\frac{p-\varepsilon}{p_k}}T^{\frac{p_k-p+\varepsilon}{p_k}}
\leq T^{\frac{p_k-p+\varepsilon}{p_k}}\leq \max\{1,T\}.
\end{equation}
On the other hand, we also have
\begin{eqnarray}\label{pqb}
\int_0^Tq\left\vert u_k\right\vert^{p-\varepsilon}\,dx&=&\int_{\{x\in[0,T]\big|\left\vert u_k(x)\right\vert\geq1\}}q\left\vert u_k\right\vert^{p-\varepsilon}\,dx+\int_{\{x\in[0,T]\big|\left\vert u_k(x)\right\vert<1\}}q\left\vert u_k\right\vert^{p-\varepsilon}\,dx\nonumber\\
&\leq&\int_{\{x\in[0,T]\big|\left\vert u_k(x)\right\vert\geq1\}}q\left\vert u_k\right\vert^{p_k}\,dx+T\max_{x\in[0,T]}q(x)\nonumber\\
&\leq&\int_0^Tq\left\vert u_k\right\vert^{p_k}\,dx+T\max_{x\in[0,T]}q(x)\nonumber\\
&\leq&1+T\max_{x\in[0,T]}q(x).
\end{eqnarray}
Clearly, (\ref{cp7}) and (\ref{pqb}) show that $\{u_k\}$ is a bounded sequence in $W_T^{1,p-\varepsilon}(0,T)$ if $k$ is large enough. Passing to a
subsequence if necessary, we can assume that $u_k \rightharpoonup u$ in $W_T^{1,p-\varepsilon}(0,T)$ and hence
that $u_k \rightarrow u$ in $C^{\alpha}[0,T]$ with $\alpha=1-1/(p-\varepsilon)$ because the embedding of $W_T^{1,p-\varepsilon}(0,T)\hookrightarrow
C^{\alpha}[0,T]$ is compact. Thus,
\begin{equation}\label{cp8}
\left\vert u_k\right\vert^{p_k}\rightarrow\vert u\vert^{p}\,\,\text{uniformly on\,\,} [0,T].
\end{equation}

We note that (\ref{cp6}) implies that
\begin{equation}\label{cp9}
\lambda_0\left(p_k\right)\int_0^T m\left\vert u_k\right\vert^{p_k}\,dx=1
\end{equation}
for all $k\in \mathbb{N}$. Thus letting $k\rightarrow+\infty$ in (\ref{cp9}) and using (\ref{cp8}), we find
\begin{equation}\label{cp10}
\liminf_{j\rightarrow+\infty}\lambda_0\left(p_j\right)\int_0^1 m\vert u\vert^{p}\,dx=1.
\end{equation}
On the other hand, since $u_k\rightharpoonup u$ in $W_T^{1,p-\varepsilon}(0,T)$, from (\ref{cp7}) we obtain that
\begin{eqnarray}
\int_0^T\left(\left\vert u'\right\vert^{p-\varepsilon}+q\vert u\vert^{p-\varepsilon}\right)\,dx&\leq&\liminf_{k\rightarrow+\infty}\int_0^T\left(\left\vert u_k'\right\vert^{p-\varepsilon}+q\left\vert u_k\right\vert^{p-\varepsilon}\right)\,dx\nonumber\\
&\leq&\lim_{k\rightarrow+\infty}\left(\left(1-\int_0^Tq
\left\vert u_k\right\vert^{p_k}\,dx\right)^{\frac{p-\varepsilon}{p_k}}T^{\frac{p_k-p+\varepsilon}{p_k}}
+\int_0^Tq\left\vert u_k\right\vert^{p-\varepsilon}\,dx\right)\nonumber\\
&=&\left(1-\int_0^Tq\vert u\vert^{p}\,dx\right)^{\frac{p-\varepsilon}{p}}T^{\frac{\varepsilon}{p}}+\int_0^Tq\vert u\vert^{p-\varepsilon}\,dx.\nonumber
\end{eqnarray}
Now, letting $\varepsilon\rightarrow 0^+$ we find
\begin{equation}\label{cp11}
\Vert u\Vert\leq1.
\end{equation}
Hence $u\in W^{1,p}(0,T)$. While we know that
$u\in W_T^{1,p-\varepsilon}(0,T)$ for each $0<\varepsilon<\varepsilon_0$. It follows that $u(0)=u(T)$. Thus, we obtain that $u\in W_T^{1,p}(0,T)$.

Finally, combining (\ref{cp10}) and (\ref{cp11}) we obtain
\begin{equation}
\liminf_{j\rightarrow+\infty}\lambda_0\left(p_j\right)\int_0^T m\vert u\vert^{p}\,dx\geq\int_0^T\left(\left\vert u'\right\vert^{p}+q\vert u\vert^p\right)\,dx.\nonumber
\end{equation}
This relationship together with the variational characterization of $\lambda_0(p)$ implies (\ref{cp4}) and hence
(\ref{cp2}). This concludes the proof of the proposition.\qed\\
\\
\textbf{Remark 3.1.} In particular, if $m\geq0$ but $m\not\equiv 0$, the results of this section are only valid for $\lambda_0^+$.
Note that some of results in this section are new even in this case.

\section{Unilateral global bifurcation}

\quad\, From now on, we use $X$ to denote the space $W_T^{1,p}(0,T)$. We start this section by considering the following auxiliary problem
\begin{equation}\label{erh}
\left\{
\begin{array}{l}
-\left(\left\vert u'\right\vert^{p-2}u'\right)'+q\varphi_p(u)=h(x),\,\, \text{ a.e.}\,\, x\in (0,T),\\
u(0)=u(T),\,\,u'(0)=u'(T)
\end{array}
\right.
\end{equation}
for a given $h\in X^*$.
\\ \\
\textbf{Lemma 4.1.} \emph{If $h\in X^*$, then problem (\ref{erh}) has a unique weak solution.}
\\ \\
\textbf{Proof.} For any $v\in X$, we define $\langle h,v\rangle:=\int_0^T hv\,dx$. It is easy to verify
that $h$ is a continuous linear functional on $X$. Since $L$ is a homeomorphism, (\ref{erh}) has a unique
solution.\qed\\

Let $G_p(h)$ denote the unique solution to problem (\ref{erh})
for a given $h\in X^*$.
Proposition 2.2 implies that $G_p:X^*\rightarrow X$ is continuous. Since $X$ embeds compactly into
$L^q(0,T)$ for each $q\in[1,+\infty]$ it follows that the restriction of $G_p$ to $L^{q'}(0,T)$ is a completely continuous operator,
where $q'=q/(q-1)$ ($+\infty$) if $q>1$ ($q=1$).
Define $T_\lambda^p(u)=G_p\left(F(\lambda,u)\right)$ on $X$, where $F(\lambda,u)$ denotes the
usual Nemitsky operator associated to $\lambda m(x)\varphi_p(u(x))$.
The compact embedding of $X\hookrightarrow L^p(0,T)$ and Theorem 1.7 of [\ref{AM}]
imply that $T_\lambda^p:X\rightarrow L^{p'}(0,T)$ is completely continuous.
Thus, $T_\lambda^p:X\rightarrow X$ is completely continuous.
Let $\Psi_{p,\lambda}$ defined on $X$ be given by
\begin{equation}
\Psi_{p,\lambda}(u)=u-T_\lambda^p(u).\nonumber
\end{equation}
Clearly, $\Psi_{p,\lambda}$ is a nonlinear
compact perturbation of the identity. Thus the Leray-Schauder degree
$\deg\left(\Psi_{p,\lambda}, B_r(0),0\right)$ is well-defined for
arbitrary $r$-ball $B_r(0)$ and $\lambda\in \left(\delta^-,\delta^+\right)\setminus\left\{\lambda_0^+(p), \lambda_0^-(p)\right\}$,
where $\delta^+$ and $\delta^-$ come from Proposition 3.1.\\

\indent Firstly, we
can compute $\deg\left(\Psi_{2,\lambda}, B_r(0),0\right)$ for any $r>0$ as the following.
\\ \\
\noindent\textbf{Lemma 4.2.} \emph{For $r>0$, we have
\begin{equation}
\deg\left(\Psi_{2,\lambda}, B_r(0),0\right)=\left\{
\begin{array}{l}
 1, \  \ \quad  ~\text{if}\ \lambda\in
\left(\lambda_0^-(2),
\lambda_0^+(2)\right),\\
-1, \ \ \   \  \text{if}\ \lambda\in \left(\lambda_0^+(2), \delta^+_2\right)\cup\left(\delta^-_2, \lambda_{0}^-(2)\right),
\end{array}
\right.\nonumber
\end{equation}
where $\delta^+_2$ and $\delta^-_2$ are chosen in such a way that there is no other eigenvalue in
$\left(\lambda_0^+(2), \delta^+_2\right)\cup\left(\delta^-_2, \lambda_{0}^-(2)\right)$ of problem (\ref{pe}) with $p=2$
.}\\
\\
\noindent\textbf{Proof.} We divide the proof into two cases.

{\it Case 1.} $\lambda\geq 0$.

Since $T_\lambda^2$ is compact and linear, by [\ref{De}, Theorem 8.10], we have
\begin{equation}
\deg\left(\Psi_{2,\lambda}, B_r(0),0\right)=(-1)^{m(\lambda)},\nonumber
\end{equation}
where $m(\lambda)$ is the sum of algebraic multiplicity of the
eigenvalues $\mu$ of problem (\ref{pe}) satisfying $\lambda^{-1}\mu<1$.
If $\lambda\in [0, \lambda^+_0(2))$, then Theorem 2 of [\ref{C1}] implies that there is no such a $\mu$ at
all, then
\begin{equation}
\deg\left(\Psi_{2,\lambda}, B_r(0),0\right)=(-1)^{m(\lambda)}=(-1)^0=1.\nonumber
\end{equation}
If $\lambda\in \left(\lambda_0^+(2), \delta^+_2\right)$, then there is only
\begin{equation}
\lambda^+_0(2)\lambda^{-1}<1.\nonumber
\end{equation}
This together with Theorem 2 of [\ref{C1}] implies
\begin{equation}
\deg\left(\Psi_{2,\lambda}, B_r(0),0\right)=-1.\nonumber
\end{equation}

{\it Case 2.} $\lambda<0$.

In this case, we consider a new sign-changing eigenvalue problem
\begin{equation}
\left\{
\begin{array}{l}
-u''+q(x)u=\widehat{\lambda} \widehat{m}(x)u,\,\,0<x<T,\\
u(0)=u(T),\,\, u'(0)=u'(T),
\end{array}
\right.\nonumber
\end{equation}
where $\hat\lambda=-\lambda$, $\hat m(x)=-m(x)$. It is easy to check
that
\begin{equation}
\hat\lambda_k^+(2)=-\lambda_k^-(2), \ \  k\in
\mathbb{N}.\nonumber
\end{equation}
Thus, we may use the result obtained in \emph{Case 1} to deduce the desired
result. \qed\\

\indent As far as the general $p$, we can compute it through the deformation along $p$.\\ \\
\noindent\textbf{Lemma 4.3.} \emph{Let $\lambda$ be a constant with
$\lambda\in\left(\delta^-,\delta^+\right)\setminus\left\{\lambda_0^+(p), \lambda_0^-(p)\right\}$. Then for
arbitrary $r>0$,}
\begin{equation}
\deg \left(I-T^p_\lambda, B_r(0),0 \right)=\left\{
\begin{array}{l}
1, \,\,\,\,\,\,\text{if}\,\,\lambda\in
\left(\lambda_0^-(p),\lambda_0^+(p)\right),\\
-1,\,\,\text{if}\,\,\lambda\in\left(\lambda_0^+(p),\delta^+\right)\cup \left(\delta^-, \lambda_0^-(p)\right).
\end{array}
\right.\nonumber
\end{equation}
\noindent\textbf{Proof.} We shall only prove for the case $\lambda>\lambda_0^+(p)$
since the proofs of other cases are completely analogous. We
also only give the proof for the case $p>2$ because proof for the case
$1<p<2$ is similar. Assume that $\lambda_{0}^+(p) < \delta^+$. Since the eigenvalue $\lambda_{0}^+(p)$
depends continuously on $p$, there exist a continuous function
$\chi:[2,p]\rightarrow\mathbb{R}$ and $q\in [2,p]$ such that
$\lambda_{0}^+(q) < \chi(q) < \delta^+$ and $\lambda=\chi(p)$.
Define
\begin{equation}
\Upsilon(q,u)=u-T_\chi^q\left(\chi(q)m(x)\varphi_q(u)\right).\nonumber
\end{equation}
It is easy to show that $\Upsilon(q,u)$ is a compact perturbation of the identity for all $u\neq 0$,
by definition of $\chi(q)$,
$\Upsilon(q,u)\neq0$, for all $q\in [2,p]$. Hence the invariance of the degree under homotopy and Lemma 4.2 imply
\begin{equation}
\deg\left(\Psi_{p,\lambda}, B_r(0),0\right)=\deg\left(\Psi_{2,\lambda},
B_r(0),0\right)=-1.\nonumber
\end{equation}
This completes the proof. \qed\\

\indent Define the Nemitskii operator $H_\lambda:X\rightarrow L^{1}(0,T)$
by
\begin{equation}
H_\lambda(u)(x):=\lambda m(x)\varphi_p(u(x))+g(x,u(x),\lambda).\nonumber
\end{equation}
Then it is clear that $H_\lambda$ is continuous operator and problem (\ref{pb}) can be equivalently written as
\begin{equation}
u=G_p\circ H_\lambda(u):=F_\lambda(u).\nonumber
\end{equation}
Since $G_p:L^1(0,T)\rightarrow X$ is compact, $F_\lambda:X\rightarrow X$ is completely continuous. Obviously, $F_\lambda(0)=0$, $\forall \lambda\in \mathbb{R}$.\\

Using a similar method to prove [\ref{DM}, Theorem 2.1] with obvious changes, we may obtain the
following result.
\\ \\
\textbf{Theorem 4.1.} \emph{$\left(\lambda_0^\nu,0\right)$ is a bifurcation
point of problem (\ref{pb}) and the corresponding bifurcation branch $\mathscr{C}_\nu$ in $\mathbb{R}\times X$
whose closure contains $\left(\lambda_0^\nu, 0\right)$ is either unbounded or contains a pair ($\overline{\lambda}, 0$)
where $\overline{\lambda}$ is an eigenvalue of problem (\ref{pe}) and $\overline{\lambda}\neq \lambda_0^\nu$.}
\\ \\
\textbf{Remark 4.1.} It is not difficult to verify that the conclusion of Lemma 2.1
is also valid for problem (\ref{pb}). It follows that
$u$ is also a classical solution for any $(\lambda,u)\in \mathscr{C}_\nu$.\\

Next, we shall prove that the first choice of the alternative of Theorem 4.1 is
the only possibility. Let $P^+$ denote the set of functions
in $X$ which are positive in $[0,T]$. Set $P^-=-P^+$ and $P =P^+\cup P^-$.
It is clear that $P^+$ and $P^-$ are disjoint and open in $X$.
Finally, let $K^{\pm}=\mathbb{R}\times P^{\pm}$ and $K=\mathbb{R}\times P$
under the product topology. \\ \\
\textbf{Lemma 4.4.} \emph{The last alternative of Theorem 4.1 is
impossible if
$\mathscr{C}_\nu\subset \left(K\cup\{\left(\lambda_0^\nu,0\right)\}\right)$.}
\\ \\
\textbf{Proof.} Suppose on the contrary, if there exists $\left(\lambda_n,u_n\right)
\rightarrow\left(\overline{\lambda},0\right)$
when $n\rightarrow+\infty$ with $\left(\lambda_n,u_n\right)\in \mathscr{C}_\nu$,
$u_n \not\equiv 0$, $\overline{\lambda}\neq \lambda_0^\nu$.
Let $v_n :=u_n/\left\Vert u_n\right\Vert$, then $v_n$ should be the solutions of the problem
\begin{equation}
v(x)=G_p\left(\lambda_n m(x)\varphi_p\left(v(x)\right)+\frac{g\left(x,u_n(x),\lambda_n\right)}{\left\Vert u_n\right\Vert^{p-1}}\right).\nonumber
\end{equation}
By an argument similar to that of [\ref{DM}, Theorem 2.1], we obtain that for some convenient
subsequence, $v_n\rightarrow v_0$ as $n\rightarrow+\infty$. It is easy to
see that $\left(\overline{\lambda}, v_0\right)$ verifies problem (\ref{pe})
and $\left\Vert v_0\right\Vert = 1$. On the other hand,
we can easily show that the bifurcation points must be eigenvalues. Thus, $\mathscr{C}_\nu$ does not join to (0,0) because
$0$ is not the eigenvalue of problem (\ref{pe}). Clearly, Proposition 2.2 implies that (0,0) is
the only solution of problem (\ref{pb}) for $\lambda= 0$. Hence, we have $\mathscr{C}_\nu\cap\left(\{0\}\times X\right)=\emptyset$.
It follows that $\overline{\lambda}\neq \lambda_0^{-\nu}$.  Theorem 3.1 follows $v_0$ must change its sign, and as a
consequence for some $n$ large enough, $u_n$ must change sign. This is a contradiction.\qed\\
\\
\textbf{Lemma 4.5.} \emph{There is a neighborhood $\mathscr{O}$ of $\left(\lambda_0^\nu,0\right)$ such that $\mathscr{C}_\nu\cap\mathscr{O}\subset \left(K\cup\{\left(\lambda_0^\nu,0\right)\}\right)$.}
\\ \\
\textbf{Proof.} If the result doesn't hold, then there would be a sequence
$\left\{\left(\lambda_n,u_n\right)\right\}\in \mathscr{C}_\nu\cap\mathscr{O}$ such that $u_n\not\equiv 0$,
$u_n\not\in K$ and $\left(\lambda_n,u_n\right)\rightarrow\left(\lambda_0^\nu,0\right)$.
Let $v_n :=u_n/\left\Vert u_n\right\Vert$, then $v_n$ should be the solutions of the problem
\begin{equation}
v(x)=G_p\left(\lambda_n m(x)\varphi_p\left(v(x)\right)+\frac{g\left(x,u_n(x),\lambda_n\right)}{\left\Vert u_n\right\Vert^{p-1}}\right).\nonumber
\end{equation}
By an argument similar to that of Lemma 4.4, we obtain for some convenient
subsequence, $v_n\rightarrow v_0$ as $n\rightarrow+\infty$. It is easy to
see that $\left(\lambda_0^\nu, v_0\right)$ verifies problem (\ref{pe})
and $\left\Vert v_0\right\Vert = 1$. Then Theorem 3.1 implies that $v_0$ is positive or negative. Without loss of generality,
we may assume that $v_0>0$ on $[0,T]$. This is impossible since $K$ is open.\qed
\\

Furthermore, applying the similar method to prove [\ref{L1}, Lemma
6.4.1] with obvious changes, we may obtain the following result, which localizes the possible solutions of (\ref{pb})
bifurcating from $\left(\lambda_0^\nu,0\right)$.\\
\\
\textbf{Lemma 4.6.} \emph{There exists a neighborhood $\mathscr{O}$ of $\left(\lambda_0^\nu,0\right)$ such that
$(\lambda,u)\in\left(\mathscr{C}_\nu\cap\mathscr{O}\right)\setminus\left\{\left(\lambda_0^\nu,0\right)\right\}$ implies $(\lambda,u)
=\left(\lambda_0^\nu+o(1), \alpha \varphi_0^\nu+y\right)$,
where $\varphi_0^\nu$ is the eigenfunction corresponding to $\lambda_0^\nu$ with $\left\Vert \varphi_0^\nu\right\Vert=1$, $\alpha\neq0$ and $y=o(\alpha)$
at $\alpha=0$.}
\\

Next, we give an important lemma which will be used later.
\\ \\
\textbf{Lemma 4.7.} \emph{If $(\lambda, u)$ is a solution of problem (\ref{pb})
and $u$ has a double zero, then $u \equiv 0$.}
\\ \\
\textbf{Proof.} Let $u$ be a solution of problem (\ref{pb}) and $x^*\in[0, T]$ be a double zero.
We note that
\begin{equation}
u(x)=\int_{x_*}^{x}\varphi_p^{-1}\left(\int_{x_*}^{s}\left(q(\tau)\varphi_p(u(\tau))-\lambda
m(\tau)\varphi_p(u(\tau))-g(\tau,u(\tau),\lambda)\right)\,d\tau\right)\,ds.\nonumber
\end{equation}
Firstly, we consider $x\in[0, x^*]$. Then we have
\begin{eqnarray}
\vert u(x)\vert&\leq&\int_x^{x^*}\varphi_p^{-1}\left(\left\vert\int_{x_*}^{s}\left(q(\tau)\varphi_p(u(\tau))
-\lambda m(\tau)\varphi_p(u(\tau))
-g(\tau,u(\tau),\lambda)\right)\,d\tau\right\vert\right)\,ds\nonumber\\
&\leq&T\varphi_p^{-1}
\left(\int_x^{x^*}
 \left\vert\left(q(\tau)\varphi_p(u(\tau))-\lambda m(\tau)\varphi_p(u(\tau))
-g(\tau,u(\tau),\lambda)\right)\right\vert\,d\tau\right).\nonumber
\end{eqnarray}
Furthermore, we get
\begin{eqnarray}
\varphi_p(\vert u(x)\vert)&\leq&T^{p-1}\int_x^{x^*}
 \left\vert\left(q(\tau)\varphi_p(u(\tau))-\lambda m(\tau)\varphi_p(u(\tau))
-g(\tau,u(\tau),\lambda)\right)\right\vert\,d\tau \nonumber\\
&\leq&T^{p-1}\int_x^{x^*}
 \left\vert\left(q(\tau)-\lambda m(\tau)
-\frac{g(\tau,u(\tau),\lambda)}{\varphi_p(u(\tau))}\right)\right\vert\varphi_p(u(\tau))\,d\tau \nonumber\\
&\leq&T^{p-1}\int_x^{x^*}
 \left(q(\tau)+\lambda \vert m(\tau)\vert
+\left\vert\frac{g(\tau,u(\tau),\lambda)}{\varphi_p(u(\tau))}\right\vert\right)\varphi_p(\vert u(\tau)\vert)\,d\tau.\nonumber
\end{eqnarray}
In view of (\ref{c1}), for any $\varepsilon>0$, there exists a constant $\delta>0$ such that
\begin{equation}
\vert g(x,s,\lambda)\vert\leq \varepsilon\varphi_p( \vert s\vert)\nonumber
\end{equation}
uniformly with respect to $x\in[0,T]$ and fixed $\lambda$ when $\vert s\vert\in[0,\delta]$.
Hence,
\begin{equation}
\varphi_p(\vert u(x)\vert)\leq T^{p-1}\int_x^{x^*}
\left(q+\vert \lambda m\vert+\varepsilon+\max_{s\in\left[\delta,\left\Vert u\right\Vert_\infty\right]}
\left\vert\frac{g(\tau,s,\lambda)}{\varphi_p(s)}\right\vert\right) \varphi_p(\vert u(\tau)\vert)\,d\tau.\nonumber
\end{equation}
By the Gronwall-Bellman inequality [\ref{Bre}], we get $u \equiv 0$ on $[0, x^*]$.
Similarly, using a modification of the Gronwall-Bellman inequality [\ref{ILL}, Lemma 2.2], we can get $u \equiv 0$ on $[x^*, T]$
and the proof is completed.\qed
\\ \\
\textbf{Theorem 4.2.} \emph{There exists an unbounded
continuum $\mathscr{C}_\nu\subseteq K\cup\{\left(\lambda_0^\nu,0\right)\}$ of solutions to problem (\ref{pb})
emanating from $\left(\lambda_0^\nu,0\right)$.}\\ \\
\textbf{Proof.}
Taking into account Theorem 4.1 and Lemma 4.4, we only need to prove that
$\mathscr{C}_\nu\subset \left(K\cup\left\{\left(\lambda_0^\nu,0\right)\right\}\right)$.
Suppose $\mathscr{C}_\nu\not\subset
\left(K\cup\left\{\left(\lambda_0^\nu,0\right)\right\}\right)$. Then Lemma 4.5 and 4.6 imply that there exists
$(\lambda, u)\in \left(\mathscr{C}_\nu\cap(\mathbb{R}\times \partial P)\right)$
such that $(\lambda, u) \neq
\left(\lambda_0^\nu, 0\right)$ and $\left(\lambda_n,u_n\right)\rightarrow(\lambda, u)$
with $\left(\lambda_n,u_n\right)\in \left(\mathscr{C}
\cap(\mathbb{R}\times P)\right)$.
The compact embedding of $X\hookrightarrow C[0,T]$ and $u_n\rightarrow u$ in $X$ imply that
$u\geq 0$ or $u\leq0$.  Without loss of generality,
we may assume that $u\geq0$ on $[0,T]$. Furthermore, $u\in \partial P$ implies that there exists a point $x_0\in[0,T]$ such that $u\left(x_0\right)=0$.
If $x_0\in(0,T)$, then Remark 4.1 and Lemma 4.7 implies $u\equiv 0$. If $x_0=0$ or $x_0=T$, then $u(T)=u(0)=0$ it implies $u'(0)\geq 0$, $u'(T)\leq0$.
Moreover, $u'(0)=u'(T)$ implies that $u'(0)=u'(T)=0$. Lemma 4.7 implies $u\equiv 0$. Let $w_n :=u_n/\left\Vert u_n\right\Vert$.
By a proof similar to that of Lemma 4.4, we can show that
there exists $w\in X$ such that $(\lambda, w)$ satisfies problem (\ref{pe}) and $\Vert w\Vert = 1$,
that is to say, $\lambda$ is an eigenvalue of problem (\ref{pe}). Therefore, $\left(\lambda_n,u_n\right)\rightarrow
\left(\lambda, 0\right)$ with $\left(\lambda_n,u_n\right)\in \mathscr{C}_\nu\cap (\mathbb{R}\times P)$.
This contradicts Lemma 4.4.\qed
\\

By an argument similar to prove [\ref{DM}, Theorem 3.2] with obvious changes, we
may obtain the following unilateral global bifurcation result.\\ \\
\textbf{Theorem 4.3.} \emph{There are two distinct unbounded sub-continua of solutions to problem (\ref{pb}),
$\mathscr{C}_\nu^+$ and $\mathscr{C}_\nu^-$, consisting of the bifurcation branch $\mathscr{C}_\nu$ and
\begin{equation}
\mathscr{C}_\nu^\sigma\subset \left(K^\sigma\cup\left\{\left(\lambda_0^\nu,0\right)\right\}\right),\nonumber
\end{equation}
where $\nu, \sigma\in\{+,-\}$.}
\\ \\
\indent Moreover, if we pose more strict assumption on $g$ as the following:
\\

(A1) $g:[0,T]\times [0,+\infty)\times \mathbb{R} \rightarrow [0,+\infty)$ is continuous.
\\

Then we have
\\ \\
\textbf{Theorem 4.4.} \emph{Besides the assumptions of Theorem 4.3, we also assume that} (A1) \emph{holds.
Then there is an unbounded continuum of solutions to problem (\ref{pb}),
$\mathscr{C}_\nu^+$ and
\begin{equation}
\mathscr{C}_\nu^+\subset \left(K^+\cup\left\{\left(\lambda_0^\nu,0\right)\right\}\right),\nonumber
\end{equation}
where $\nu\in\{+,-\}$.}
\\ \\
\textbf{Proof.} Define
\begin{equation}
\widetilde{g}(x,s,\lambda)=\left\{
\begin{array}{l}
g(x,s,\lambda),\,\,\,\,~~~\,\text{if\,\,}s>0,\\
0,\,\,\,\,~~~~~~~~~~~~~\,\text{if\,\,}s=0,\\
-g(x,-s,\lambda),\,\, \text{if\,\,}s<0.
\end{array}
\right.\nonumber
\end{equation}
We consider the following problem
\begin{equation}\label{ppb}
\left\{
\begin{array}{l}
-\left(\varphi_p\left(u'\right)\right)'+q(x)\varphi_p(u)=\lambda m(x)\varphi_p(u)+\widetilde{g}(x,u,\lambda),\,\,0<x<T,\\
u(0)=u(T),\,\, u'(0)=u'(T),
\end{array}
\right.
\end{equation}
Applying Theorem 4.3 to problem (\ref{ppb}), we obtain that there are two distinct unbounded sub-continua of solutions to problem (\ref{ppb}),
$\mathscr{C}_\nu^+$ and $\mathscr{C}_\nu^-$, consisting of the bifurcation branch $\mathscr{C}_\nu$ and
\begin{equation}
\mathscr{C}_\nu^\sigma\subset \left(K^\sigma\cup\{\left(\lambda_0^\nu,0\right)\}\right),\nonumber
\end{equation}
where $\nu,\sigma\in\{+,-\}$. Clearly, $\mathscr{C}_\nu^+$  is also the solution branch of problem (\ref{pb}).\qed\\
\\
\textbf{Remark 4.2.} Note that if $m\geq0$ but $m\not\equiv 0$, we can only get the component $\mathscr{C}_+^\sigma$
emanating from $\left(\lambda_0^+,0\right)$.
Thus, our results in this section are new even in the definite weight case.

\section{One-sign solutions with signum condition}

\quad\, In this section, we shall investigate the existence and multiplicity of one-sign
solutions to problem (\ref{pn}).

\indent Let $f_0, f_\infty\in \mathbb{R}\setminus \mathbb{R}^-$ such that
\begin{equation}
f_0=\lim\limits_{|s|\rightarrow 0}\frac{f(s)}{\varphi_p(s)}\,\,\text{and}\,\,f_{\infty}=\lim\limits_{|s|\rightarrow
+\infty}\frac{f(s)}{\varphi_p(s)}.\nonumber
\end{equation}
Through out this section, we always suppose that $f$ satisfies the following signum condition\\

\emph{(A2) $f\in C(\mathbb{R},\mathbb{R})$ with
$f(s)s>0$ for $s\neq0$.}\\

Clearly, (A2) implies $f(0)=0$. Hence, $u=0$ is always the solution of problem (\ref{pn}). Applying
Theorem 4.3, we shall establish the existence of one-sign
solutions of problem (\ref{pn}) as the following.\\ \\
\textbf{Theorem 5.1.} \emph{If $f_0\in(0,+\infty)$ and $f_\infty\in(0,+\infty)$,
then for any $\lambda\in\left(\lambda_0^\nu/f_\infty,\lambda_0^\nu/f_0\right)\cup\left(\lambda_0^\nu/f_0,\lambda_0^\nu/f_\infty\right)$, problem (\ref{pn})
has two solutions $u^+$ and $u^-$ such
that $u^+$ is positive and $u^-$ is negative on $[0,T]$.}
\\ \\
\indent In order to prove Theorem 5.1, we need the following Sturm-type comparison result.\\
\\
\textbf{Lemma 5.1.} \emph{Let $b_2(x)\geq b_1(x)$ for $x\in(0,T)$ and $b_i(x)\in C(0,T)$, $i=1,2$. Also let $u_1$, $u_2$
be solutions of the following differential equations
\begin{equation}
\left(\varphi_p\left(u'\right)\right)'+b_i(x)\varphi_p(u)-q(x)\varphi_p(u)=0,\,\, i=1,2,\,\, x\in(0,T),\nonumber
\end{equation}
respectively. If $(c,d)\subset(0,T)$, and $u_1(c)=u_1(d)=0$, $u_1(x)\neq 0$ in $(c,d)$, then either there exists $\tau\in (c,d)$
such that $u_2(\tau)=0$ or $b_2=b_1$ and $u_2(x)=\mu u_1(x)$ for some constant $\mu\neq0$.}
\\ \\
\textbf{Proof.} Applying a similar method to prove [\ref{DMX}, Lemma 3.1] with obvious changes, we can show this lemma.
\qed\\
\\
\indent Denote
\begin{equation}
I^+:=\left\{x\in [0,T]\,|\, m(x)>0\right\}, \ \ \ I^-:=\left\{x\in [0,T]\,|\, m(x)<0\right\}.\nonumber
\end{equation}
\noindent\textbf{Lemma 5.2.} \emph{Let $\widehat{I}=[a,b]$ such
that $\widehat{I}\subset I^+$ and $\text{meas}\,\left\{\widehat{I}\right\}>0$.
And assume that $g_n:[0,T]\to (0, +\infty)$ is a continuous function such that
\begin{equation}
\lim_{n\to +\infty} g_n(x)=+\infty\ \ \text{uniformly on}\
\widehat{I}.\nonumber
\end{equation}
Let $y_n\in X$ be a solution of the equation
\begin{equation}
\left(\varphi_p\left(y_n'\right)\right)'-q(x)\varphi_p\left(y_n\right)+ m(x)g_n(x)\varphi_p\left(y_n\right)=0, \,\, x\in (0, T).\nonumber
\end{equation}
Then the number of zeros of $y_n|_{\widehat{I}}$ goes to infinity as $n\to +\infty$.}
\\ \\
\noindent\textbf{Proof.} Taking a subsequence if
necessary, we may assume that
\begin{equation}
m(x)g_{n_j}(x)\geq \lambda_j, \ \  x\in \widehat{I}\nonumber
\end{equation}
as $j\to +\infty$, where $\lambda_j$ is the $j$th eigenvalue of the following problem
\begin{equation}\label{ppp}
\left\{
\begin{array}{l}
\left(\varphi_p\left(u'(x)\right)\right)'-q(x)\varphi_p(u(x))+\lambda \varphi_p(u(x))=0,\,\,x\in (0,T),\\
u(0)=u(T)=0.
\end{array}
\right.
\end{equation}
Set $t=(b-a)/Tx+a$, $v(t)=u\left(T/(b-a)(t-a)\right)$ and $q(t)=q\left(T/(b-a)(t-a)\right)$.
By some simple computations, we can show
\begin{equation}\label{ppp}
\left\{
\begin{array}{l}
\left(\varphi_p\left(v'(t)\right)\right)'-q(t)\varphi_p(v(t))+\lambda \varphi_p(v(t))=0,\,\,t\in (a,b),\\
v(a)=v(b)=0.
\end{array}
\right.
\end{equation}
Let $\varphi_j$ be the corresponding eigenvalue of $\lambda_j$.
Theorem 3.1 of [\ref{R1}] implies that the number of zeros of $\varphi_j\big|_{\widehat{I}}$ goes to infinity as $j\to +\infty$.
By Lemma 5.1, one obtains that the number of zeros of $y_n|_{\widehat{I}}$ goes to infinity as $n\to +\infty$.
It follows the desired results.\qed
\\ \\
\textbf{Proof of Theorem 5.1.} We only prove the case of $\lambda>0$ since the proof of $\lambda<0$ can be given similarly.
Let $\zeta\in C(\mathbb{R})$ such that $f(s)=f_0 \varphi_p(s)+\zeta(s)$
with $\lim_{s\rightarrow0}\zeta(s)/\varphi_p(s)=0.$
By Theorem 4.3, we have that
there are two distinct unbounded sub-continua $\mathscr{C}_+^+$ and $\mathscr{C}_+^-$,
consisting of the bifurcation branch $\mathscr{C}_+$ emanating
from $\left(\lambda_0^+/f_0, 0\right)$, such that
\begin{equation}
\mathscr{C}_+^\sigma\subset \left(\left\{\left(\lambda_0^+,0\right)\right\}
\cup\left(\mathbb{R}\times P^\sigma\right)\right).\nonumber
\end{equation}
To complete the proof of this theorem, it will be enough to show that $\mathscr{C}_+^\sigma$ joins
$\left(\lambda_0^+/f_0, 0\right)$ to
$\left(\lambda_0^+/f_\infty, +\infty\right)$. Let
$\left(\lambda_n, u_n\right) \in \mathscr{C}_+^\sigma$ satisfying
$\lambda_n+\left\Vert u_n\right\Vert\rightarrow+\infty.$
We note that $\lambda_n >0$ for all $n \in \mathbb{N}$ since (0,0) is
the only solution of problem (\ref{pn}) for $\lambda= 0$ and
$\mathscr{C}_+^\sigma\cap\left(\{0\}\times X\right)=\emptyset$.

We divide the rest proofs into two steps.

\emph{Step 1.} We show that there exists a constant $M$ such that $\lambda_n
\in(0,M]$ for $n\in \mathbb{N}$ large enough.

On the contrary, we suppose that $\lim_{n\rightarrow +\infty}\lambda_n=+\infty.$
We note that
\begin{equation}
-\left(\left(\varphi_p\left(u_n'\right)\right)\right)'+q\varphi_p\left(u_n\right)=\lambda_n m
\widetilde{f}_n(x)\varphi_p\left(u_n\right),\nonumber
\end{equation}
where
\begin{equation}
\widetilde{f}_n(x)=\left\{
\begin{array}{l}
\frac{f\left(u_n\right)}{\varphi_p\left(u_n\right)},\,\, \text{if}\,\,u_n\neq0,\\
f_0,\,\,\,\,\,\,\quad\text{if}\,\,u_n=0.
\end{array}
\right.\nonumber
\end{equation}
The signum condition (A2) implies that there exists a positive
constant $\varrho$ such that $\widetilde{f}_n(x)\geq\varrho$ for
any $x\in[0,T]$.
By Lemma 5.2, we get that $u_n$ must change its sign in $[0,T]$ for $n$
large enough, and this contradicts the fact that $u_n\in \mathscr{C}_+^\sigma$.

\emph{Step 2.} We show that $\mathscr{C}_+^\sigma$ joins
$\left(\lambda_0^+/f_0, 0\right)$ to
$\left(\lambda_0^+/ f_\infty, +\infty\right)$.

It follows from \emph{Step 1} that $\left\Vert u_n\right\Vert\rightarrow+\infty.$
Let $\xi\in C(\mathbb{R})$ be such that $f(s)=f_\infty\varphi_p(s)+\xi(s).$
Then $\lim_{\vert s\vert\rightarrow+\infty}\xi(s)/\varphi_p(s)=0.$
Let $\widetilde{\xi}(u)=\max_{u\leq \vert s\vert\leq 2u}\vert \xi(s)\vert.$
Then $\widetilde{\xi}$ is nondecreasing and
\begin{equation}\label{eu0}
\lim_{u\rightarrow +\infty}\frac{\widetilde{\xi}(u)}{\varphi_p(u)}=0.
\end{equation}

We divide the equation
\begin{equation}
-\left(\varphi_p\left(u_n'\right)\right)'+q\varphi_p\left(u_n\right)=\lambda_n f_\infty  \varphi_p\left(u_n\right)+\lambda_nm\xi\left(u_n\right)\nonumber
\end{equation}
by $\left\Vert u_n\right\Vert$ and set $\overline{u}_n = u_n/\left\Vert u_n\right\Vert$.
Since $\overline{u}_n$ are bounded in $X$,
after taking a subsequence if
necessary, we have that $\overline{u}_n \rightharpoonup \overline{u}$
for some $\overline{u} \in X$. Moreover, from
(\ref{eu0}) and the fact that $\widetilde{\xi}$ is nondecreasing,
we have that
\begin{equation}\label{000}
\lim_{n\rightarrow+\infty}\frac{ \xi\left(u_n\right)}{\left\Vert u_n\right\Vert^{p-1}}=0,
\end{equation}
since
\begin{equation}
\frac{ \vert\xi\left(u_n\right)\vert}{\left\Vert u_n\right\Vert^{p-1}}\leq\frac{ \widetilde{\xi}
(\left\vert u_n\right\vert)}{\left\Vert u_n\right\Vert^{p-1}}
\leq\frac{ \widetilde{\xi}(\left\Vert u_n\right\Vert_\infty)}{\left\Vert u_n\right\Vert^{p-1}}\leq\frac{ \widetilde{\xi}(C_0\left\Vert u_n\right\Vert)}{\left\Vert u_n\right\Vert^{p-1}}\leq\frac{ C_0^{p-1}\widetilde{\xi}(\left\Vert C_0 u_n\right\Vert)}{\left\Vert C_0 u_n\right\Vert^{p-1}},\nonumber
\end{equation}
where $\Vert \cdot\Vert_\infty$ denotes the usual norm of $C[0,T]$ and $C_0$ is the embedding constant of $X\hookrightarrow C[0,T]$.

By the continuity and compactness of $F_\lambda$, it follows that
\begin{equation}
-\left(\varphi_p\left(\overline{u}'\right)\right)'+q\varphi_p(\overline{u})=\overline{\lambda} f_\infty m  \varphi_p(\overline{u}),\nonumber
\end{equation}
where
$\overline{\lambda}=\underset{n\rightarrow+\infty}\lim\lambda_n$, again
choosing a subsequence and relabeling it if necessary.

It is clear that $\Vert \overline{u}\Vert=1$ and $\overline{u}\in
\overline{\mathscr{C}_+^\sigma}\subseteq
\mathscr{C}_+^\sigma$ since $\mathscr{C}_+^\sigma$ is closed in $\mathbb{R}\times X$.
Thus, $\overline{\lambda}f_\infty=\lambda_0^+$, i.e., $\overline{\lambda}=\lambda_0^+/ f_\infty$.
Therefore, $\mathscr{C}_+^\sigma$ joins $\left(\lambda_0^+/f_0, 0\right)$ to $\left(\lambda_0^+/f_\infty,+\infty\right)$.\qed\\

From the proof of Theorem 5.1, we can easily get the following corollary.
\\ \\
\textbf{Corollary 5.1.} \emph{Assume that there exists a positive constant $\rho>0$ such that
\begin{equation}
\frac{f(s)}{\varphi_p(s)}\geq \rho\nonumber
\end{equation}
for any $s\neq 0$.
Then there exist $\lambda_*^+>0$ and $\lambda_*^-<0$ such that problem (\ref{pn}) has no one-sign solution for
any $\lambda\in\left(-\infty,\lambda_*^-\right)\cup\left(\lambda_*^+,+\infty\right)$.}
\\
\\
\textbf{Theorem 5.2.} \emph{If $f_0\in(0,+\infty)$ and $f_\infty=0$,
then for any $\lambda\in\left(\lambda_0^+/f_0,+\infty\right)\cup\left(-\infty,\lambda_0^-/f_0\right)$, problem (\ref{pn})
has two solutions $u^+$ and $u^-$ such
that $u^+$ is positive and $u^-$ is negative on $[0,T]$.}
\\ \\
\textbf{Proof.} We shall only prove the case $\lambda>0$
since the proof for the other case is completely analogous. In view of Theorem 5.1, we only need to show that
$\mathscr{C}_+^\sigma$ joins
$\left(\lambda_0^+/f_0, 0\right)$ to
$\left(+\infty, +\infty\right)$.
Suppose on the contrary that there exists $\lambda_M$ be a blow up point and $\lambda_M<+\infty$.
Then there exists a sequence $\left\{\lambda_n, u_n\right\}$ such that
$\underset{n\rightarrow +\infty}{\lim}\lambda_n=\lambda_{M}$ and
$\underset{n\rightarrow +\infty}{\lim}\left\Vert u_n\right\Vert=+\infty$ as
$n\rightarrow+\infty$. Let $v_n =u_n/\left\Vert u_n\right\Vert$. Then $v_n$
should be the solutions of problem
\begin{equation}
v=G_p\left(\frac{\mu_n m(x)f\left(u_n(x)\right)}{\left\Vert u_n\right\Vert^{p-1}}\right).\nonumber
\end{equation}
Similar to (\ref{000}), we can show
\begin{equation}
\lim_{n\rightarrow+\infty}\frac{ f\left(u_n\right)}{\left\Vert u_n\right\Vert^{p-1}}=0.\nonumber
\end{equation}
By the compactness of $F_\lambda$, we obtain that for some convenient subsequence
$v_n\rightarrow v_0$ as $n\rightarrow+\infty$. Letting $n\rightarrow+\infty$,
we obtain that
$v_0\equiv0$. This contradicts $\left\Vert v_0\right\Vert=1$. \qed\\

\indent Next, we shall need the following topological lemma.
\\ \\
\textbf{Lemma 5.3 (see [\ref{MA}].} \emph{Let $X$ be a Banach space and let $C_n$ be a family of closed connected subsets of $X$. Assume that:}
\\

(i) \emph{there exist $z_n\in C_n$, $n=1,2,\ldots$, and $z^*\in X$, such that $z_n\rightarrow z^*$;}

(ii) \emph{$r_n=\sup \left\{\Vert x\Vert\big| x\in C_n\right\}=+\infty$;}

(iii) \emph{for every $R>0$, $\left(\cup_{n=1}^{+\infty} C_n\right)\cap B_R$ is a relatively compact set of $X$, where}
\begin{equation}
B_R=\{x\in X|\Vert x\Vert\leq R\}.\nonumber
\end{equation}

\noindent \emph{Then there exists an unbounded component $\mathfrak{C}$ of $\mathfrak{D} =: \limsup_{n\rightarrow +\infty}C_n$ and $z\in \mathfrak{C}$.}
\\ \\
\textbf{Theorem 5.3.} \emph{If $f_0\in(0,+\infty)$ and $f_\infty=+\infty$,
then for any $\lambda\in\left(0,\lambda_0^+/f_0\right)\cup\left(\lambda_0^-/f_0,0\right)$, problem (\ref{pn}) has two
solutions $u^+$ and $u^-$ such
that $u^+$ is positive and $u^-$ is negative on $[0,T]$.}
\\ \\
\textbf{Proof.} Inspired by the idea of [\ref{ACD}], we define the cut-off function of $f$ as the following
\begin{equation}
f_n(s)=\left\{
\begin{array}{l}
f(s),\,\,\quad\quad\quad\quad\quad\quad\,\,
\quad\quad\quad\,\,\quad s\in\left[-n,n\right],\\
\frac{n\varphi_p(2n)-f(n)}{n}(s-n)+f(n),\quad\,\,\,\,s\in\left(n,2n\right),\\
\frac{n\varphi_p(2n)+f(-n)}{n}(s+n)+f(-n),\,\,s\in\left(-2n,-n\right),\\
n\varphi_p(s),\quad\quad\quad\quad\quad\quad\quad\quad\,\,
\quad\,\,\,\,s\in\left(-\infty,-2n\right]\cup\left[2n,+\infty\right).
\end{array}
\right.\nonumber
\end{equation}
We consider the following problem
\begin{equation}\label{pn3}
\left\{
\begin{array}{l}
-\left(\varphi_p\left(u'\right)\right)'+q(x)\varphi_p(u)=\lambda m(x)f_n(u),\,\,0<x<T,\\
u(0)=u(T),\,\, u'(0)=u'(T).
\end{array}
\right.
\end{equation}
Clearly, we can see that $\lim_{n\rightarrow+\infty}f_n(s)=f(s)$, $\left(f_n\right)_0=f_0$ and $\left(f_n\right)_\infty=n$.
Theorem 5.1 implies
that there exists a sequence of unbounded continua $\left(\mathscr{C}_\nu^\sigma\right)_n$ of solutions to problem (\ref{pn3})
emanating from $\left(\lambda_0^\nu/f_0, 0\right)$
and joining to $\left(\lambda_0^\nu/n, +\infty\right)$.

By Lemma 5.3, there exists an unbounded component $\mathscr{C}_\nu^\sigma$ of $\lim\sup_{n\rightarrow+\infty} \left(\mathscr{C}_\nu^\sigma\right)_n$
such that $\left(\lambda_0^\nu/f_0, 0\right)\in\mathscr{C}_\nu^\sigma$ and $\left(0, +\infty\right)\in\mathscr{C}_\nu^\sigma$.
This completes the proof.\qed
\\ \\
\textbf{Theorem 5.4.} \emph{If $f_0=0$ and $f_\infty\in(0,+\infty)$,
then for any $\lambda\in\left(\lambda_0^+/f_\infty,+\infty\right)\cup \left(-\infty,\lambda_0^-/f_\infty\right)$, problem (\ref{pn})
has two solutions $u^+$ and $u^-$ such
that $u^+$ is positive and $u^-$ is negative on $[0,T]$.}
\\ \\
\textbf{Proof.} If $(\lambda,u)$ is any nontrivial solution of problem (\ref{pn}), dividing problem (\ref{pn}) by $\Vert u\Vert^{2(p-1)}$ and
setting $v=u/\Vert u\Vert^2$ yields
\begin{equation}\label{pn4}
\left\{
\begin{array}{l}
-\left(\varphi_p\left(v'\right)\right)'+q(x)\varphi_p(v)=\lambda m(x)\frac{f(u)}{\Vert u\Vert^{2(p-1)}},\,\,0<x<T,\\
v(0)=v(T),\,\, v'(0)=v'(T).
\end{array}
\right.
\end{equation}
Define
\begin{equation}
\widetilde{f}(v)=\left\{
\begin{array}{l}
\Vert v\Vert^{2(p-1)}f\left(\frac{v}{\Vert v\Vert^2}\right),\,\,\text{if}\,\, v\neq 0,\\
0,\,\ \quad\quad\quad\quad\quad\quad\quad\text{if}\,\, v=0.
\end{array}
\right.\nonumber
\end{equation}

Evidently, problem (\ref{pn4}) is equivalent to
\begin{equation}\label{pn44}
\left\{
\begin{array}{l}
-\left(\varphi_p\left(v'\right)\right)'+q(x)\varphi_p(v)=\lambda m(x)\widetilde{f}(v),\,\,0<x<T,\\
v(0)=v(T),\,\, v'(0)=v'(T).
\end{array}
\right.
\end{equation}
It is obvious that $(\lambda,0)$ is always the solution of problem (\ref{pn44}).
By simple computation, we can show that $\widetilde{f}_0=f_\infty$
and $\widetilde{f}_\infty=f_0$.

Now, applying Theorem 5.2 and the inversion $v\rightarrow v/\Vert v\Vert^2=u$,
we achieve the conclusion.\qed\\
\\
\textbf{Theorem 5.5.} \emph{If $f_0=+\infty$ and $f_\infty\in(0,+\infty)$,
then for any $\lambda\in\left(0,\lambda_0^+/f_\infty\right)\cup\left(\lambda_0^-/f_\infty,0\right)$, problem (\ref{pn}) has
two solutions $u^+$ and $u^-$ such that $u^+$ is positive and $u^-$ is negative on $[0,T]$.}
\\ \\
\textbf{Proof.} By an argument similar to that of Theorem 5.4 and the conclusion of Theorem 5.3, we can obtain the conclusion.\qed
\\
\\
\textbf{Theorem 5.6.} \emph{If $f_0=0$ and $f_\infty=+\infty$,
then for any $\lambda\in\left(0,+\infty\right)\cup\left(-\infty,0\right)$, problem (\ref{pn}) has two solutions $u^+$ and $u^-$ such
that $u^+$ is positive and $u^-$ is negative on $[0,T]$.}
\\ \\
\textbf{Proof.} Define
\begin{equation}
f^n(s)=\left\{
\begin{array}{l}
\frac{1}{n}\varphi_p(s),\,\,\quad\quad\quad\quad\quad\quad\quad\,\,\,\,\,\,\,\,
\quad\quad\quad\,\,\quad s\in\left[-\frac{1}{n},\frac{1}{n}\right],\\
\left(f\left(\frac{2}{n}\right)-\frac{1}{n^p}\right)(ns-2)+f\left(\frac{2}{n}\right),\quad\quad\,\,\,\,\,s\in\left(\frac{1}{n},\frac{2}{n}\right),\\
-\left(f\left(-\frac{2}{n}\right)+\frac{1}{n^{p}}\right)(ns+2)+f\left(-\frac{2}{n}\right),\,\,s\in\left(-\frac{2}{n},-\frac{1}{n}\right),\\
f(s),\quad\quad\quad\quad\quad\quad\quad\quad\quad\quad\quad\quad\,\,
\quad\,\,\,\,s\in\left(-\infty,-\frac{2}{n}\right]\cup\left[\frac{2}{n},+\infty\right).
\end{array}
\right.\nonumber
\end{equation}
Now, consider the following problem
\begin{equation}\label{pn6}
\left\{
\begin{array}{l}
-\left(\varphi_p\left(u'\right)\right)'+q(x)\varphi_p(u)=\lambda m(x)f^n(u),\,\,0<x<T,\\
u(0)=u(T),\,\, u'(0)=u'(T).
\end{array}
\right.
\end{equation}
It is obviou that $\lim_{n\rightarrow+\infty}f^n(s)=f(s)$, $f_0^n=1/n$ and $f_\infty^n=f_\infty=\infty$.
Theorem 5.3 implies
that there exists a sequence of unbounded components $\left(\mathscr{C}_\nu^\sigma\right)_n$ of solutions to problem (\ref{pn6})
emanating from $\left(n\lambda_0^\nu, 0\right)$
and joining to $\left(0, +\infty\right)$.

Lemma 5.3 implies that there exists an unbounded component $\mathscr{C}_\nu^\sigma$ of $\lim\sup_{n\rightarrow+\infty} \left(\mathscr{C}_\nu^\sigma\right)_n$
such that $(0,+\infty)\in\mathscr{C}_\nu^\sigma$ and $\left(+\infty, 0\right)\in\mathscr{C}_\nu^\sigma$.\qed
\\ \\
\textbf{Theorem 5.7.} \emph{If $f_0=+\infty$ and $f_\infty=0$,
then for any $\lambda\in\left(0,+\infty\right)\cup (-\infty,0)$, problem (\ref{pn}) has two solutions $u^+$ and $u^-$ such
that $u^+$ is positive and $u^-$ is negative on $[0,T]$.}
\\ \\
\textbf{Proof.} By an argument similar to that of Theorem 5.4 and the conclusions of Theorem 5.6, we can prove the conclusion.\qed
\\ \\
\textbf{Theorem 5.8.} \emph{If $f_0=0$ and $f_\infty=0$,
then there exist $\lambda_+^+>0$ and $\lambda_-^+<0$ such that for any $\lambda\in\left(\lambda_+^+,+\infty\right)\cup\left(-\infty,\lambda_-^+\right)$,
problem (\ref{pn}) has two positive solutions $u_1^+$ and $u_2^+$ on $[0,T]$. Similarly, there exist $\lambda_+^->0$ and $\lambda_-^-<0$
such that for any $\lambda\in\left(\lambda_+^-,+\infty\right)\cup\left(-\infty,\lambda_-^-\right)$, problem (\ref{pn}) has two
negative solutions $u_1^-$ and $u_2^-$ on $[0,T]$. Moreover, there exists $\mu_*^\nu>0$ such that problem (\ref{pn}) has no one-sign solutions for any $\lambda\in\left(0,\mu_*^\nu\right)$.}
\\ \\
\textbf{Proof.} Define
\begin{equation}
g_n(s)=\left\{
\begin{array}{l}
\frac{1}{n}\varphi_p(s),\,\,\quad\quad\quad\quad\quad\quad\quad\,\,\,\,\,\,\,\,
\quad\quad\quad\,\,\quad s\in\left[-\frac{1}{n},\frac{1}{n}\right],\\
\left(f\left(\frac{2}{n}\right)-\frac{1}{n^p}\right)(ns-2)+f\left(\frac{2}{n}\right),\quad\quad\,\,\,\,\,s\in\left(\frac{1}{n},\frac{2}{n}\right),\\
-\left(f\left(-\frac{2}{n}\right)+\frac{1}{n^{p}}\right)(ns+2)+f\left(-\frac{2}{n}\right),\,\,s\in\left(-\frac{2}{n},-\frac{1}{n}\right),\\
f(s),\quad\quad\quad\quad\quad\quad\quad\quad\quad\quad\quad\quad\,\,
\quad\,\,\,\,s\in\left(-\infty,-\frac{2}{n}\right]\cup\left[\frac{2}{n},+\infty\right).
\end{array}
\right.\nonumber
\end{equation}
By the conclusions of Theorem 5.2 and an argument similar to that of Theorem 5.6, we can obtain an unbounded component $\mathscr{C}_\nu^\sigma$ of solutions to
problem (\ref{pn}) such that $(+\infty, 0)\in\mathscr{C}_\nu^\sigma$ and $(+\infty,+\infty)\in\mathscr{C}_\nu^\sigma$.

Finally, we show that there exists $\mu_*^\nu>0$ such that problem (\ref{pn}) has no one-sign solutions for any $\lambda\in\left(0,\mu_*^\nu\right)$.
Suppose on the contrary that there exists a sequence $\left\{\lambda_n, u_n\right\}\in\mathscr{C}_\nu^\sigma$ such that
$\underset{n\rightarrow +\infty}{\lim}\mu_n=0$. $f_0=f_\infty=0$ implies that
there exists a positive constant $M$ such that
\begin{equation}
\left\vert\frac{f(s)}{\varphi_p(s)}\right\vert\leq M \,\,\text{for any}\,\, s\neq 0.\nonumber
\end{equation}
Let $v_n =u_n/\left\Vert u_n\right\Vert$. Obviously, one has
\begin{equation}
v_n=G_p\left(\frac{\lambda_n m(x)f\left(u_n(x)\right)}{\left\Vert u_n\right\Vert^{p-1}}\right).\nonumber
\end{equation}
By the compactness of $F_\lambda$, we obtain that for some convenient subsequence
$v_n\rightarrow v_1$ as $n\rightarrow+\infty$. Letting $n\rightarrow+\infty$,
we obtain that $v_0\equiv0$. This contradicts $\left\Vert v_1\right\Vert=1$.
\qed\\

From the proof of Theorem 5.8, we can deduce the following corollary.
\\ \\
\textbf{Corollary 5.2.} \emph{Assume that there exists a positive constant $\varrho>0$ such that
\begin{equation}
\left\vert\frac{f(s)}{\varphi_p(s)}\right\vert\leq\varrho\nonumber
\end{equation}
for any $s\neq 0$.
Then there exist $\mu_*^+>0$ and $\mu_*^-<0$ such that problem (\ref{pn}) has no one-sign solution for
any $\lambda\in\left(0,\mu_*^-\right)\cup\left(0,\mu_*^+\right)$.}
\\ \\
\textbf{Theorem 5.9.} \emph{If $f_0=+\infty$ and $f_\infty=+\infty$,
then there exist $\mathcal{\lambda}_+^+>0$ and $\mathcal{\lambda}_-^+<0$ such that for any
$\lambda\in\left(0,\lambda_+^+\right)\cup\left(\lambda_-^+,0\right)$, problem (\ref{pn}) has two positive solutions
$u_1^+$ and $u_2^+$ on $[0,T]$. Similarly, there exist $\lambda_+^->0$
and $\lambda_-^-<0$ such that for any $\lambda\in\left(0,\lambda_+^-\right)\cup\left(\lambda_-^-,0\right)$, problem (\ref{pn})
has two negative solutions $u_1^-$
and $u_2^-$ on $[0,T]$.}\\
\\
\textbf{Proof.} Define
\begin{equation}
g^n(s)=\left\{
\begin{array}{l}
n\varphi_p(s),\,\,\quad\quad\quad\quad\quad\quad\quad\,\,\,~\,\,
\quad\quad\quad\quad\,\,\quad s\in\left[-\frac{1}{n},\frac{1}{n}\right],\\
\left(f\left(\frac{2}{n}\right)-\frac{1}{n^{p-2}}\right)(ns-2)+f\left(\frac{2}{n}\right),\quad\quad\,\,\,\,\,s\in
\left(\frac{1}{n},\frac{2}{n}\right),\\
-\left(f\left(-\frac{2}{n}\right)+\frac{1}{n^{p-2}}\right)(ns+2)+f\left(-\frac{2}{n}\right),\,\,s\in
\left(-\frac{2}{n},-\frac{1}{n}\right),\\
f(s),\quad\quad\quad\quad\quad\quad\quad\quad\quad\,
\quad\quad\quad\quad\quad\,\,\,\,s\in\left(-\infty,-\frac{2}{n}\right]\cup\left[\frac{2}{n},+\infty\right).
\end{array}
\right.\nonumber
\end{equation}
We consider the following problem
\begin{equation}\label{pn9}
\left\{
\begin{array}{l}
-\left(\varphi_p\left(u'\right)\right)'+q(x)\varphi_p(u)=\lambda m(x)g^n(u),\,\,0<x<T,\\
u(0)=u(T),\,\, u'(0)=u'(T).
\end{array}
\right.
\end{equation}
It is no difficulty to verify that $\lim_{n\rightarrow+\infty}g^n(s)=f(s)$, $g_0^n=n$ and $g_\infty^n=f_\infty=\infty$.
Theorem 5.3 implies
that there exists a sequence of unbounded continua $\left(\mathscr{C}_\nu^\sigma\right)_n$ of solutions to problem (\ref{pn9})
emanating from $\left(\lambda_0^\nu/n, 0\right)$
and joining to $\left(0, +\infty\right)$.

By making use of Lemma 5.3 again, we obtain that there exists an unbounded component $\mathscr{C}_\nu^\sigma$ of $\lim\sup_{n\rightarrow+\infty} \left(\mathscr{C}_\nu^\sigma\right)_n$
such that $(0,+\infty)\in\mathscr{C}_\nu^\sigma$ and $\left(0, 0\right)\in\mathscr{C}_\nu^\sigma$.\qed\\

Now, we strengthen the assumptions on $f$ and $m$ as the following\\

(A3) $f:[0,+\infty)\rightarrow[0,+\infty)$ is continuous and $f(s)>0$ for $s>0$;

(A4) $m:[0,T]\rightarrow[0,+\infty)$ is continuous and $m\not\equiv 0$.\\

By Theorem 4.4, Remark 4.2, Theorem 5.1--5.9 and Corollary 5.1--5.2, we can easily show the following corollary.
\\ \\
\textbf{Corollary 5.3.} \emph{Assume that} (A3)--(A4) \emph{hold}.\\

(a) \emph{If $f_0=0$ or $f_\infty=0$, then there exists $\lambda_0>0$ such that (\ref{pn}) has a positive solution for $\lambda>\lambda_0$.}

(b) \emph{If $f_0=+\infty$ or $f_\infty=+\infty$, then there exists $\lambda_0>0$ such that (\ref{pn}) has a positive solution for $0<\lambda<\lambda_0$.}

(c) \emph{If $f_0=f_\infty=0$, then there exists $\lambda_0>0$ such that (\ref{pn}) has at least two positive solutions for $\lambda>\lambda_0$.}

(d) \emph{If $f_0=f_\infty=+\infty$, then there exists $\lambda_0>0$ such that (\ref{pn}) has at least two positive solutions for $0<\lambda<\lambda_0$.}

(e) \emph{If $f_0<+\infty$ and $f_\infty<+\infty$, then there exists $\lambda_0>0$ such that (\ref{pn}) has no positive solutions for $0<\lambda<\lambda_0$.}

(f) \emph{If $f_0>0$ and $f_\infty>0$, then there exists $\lambda_0>0$ such that (\ref{pn}) has no positive solutions for $\lambda>\lambda_0$.}
\\ \\
%
%
%
%
\textbf{Remark 5.1.} Note that the solutions obtained from Theorem 5.1--5.9 are also classical solutions by Remark 4.1.
\\ \\
\textbf{Remark 5.2.} We also note that if $f$ is singular at 0 then $f_0=+\infty$. Thus, singular nonlinearity is a special
case in Theorem 5.5, 5.7 and 5.9.

\section{Uniqueness of positive solutions}

\quad\, In this section, under some more strict assumptions of $f$, we shall
show that the unbounded continua which are obtained in Section 5
may be curves. We just show the case of $f_0=+\infty$ and $f_\infty=0$. Other cases can be discussed similarly.\\

Firstly, we give the following assumption:
\\

(A5) $f(s)/\varphi_p(s)$ is strictly decreasing in $(0,+\infty)$.\\

Under the assumptions (A3)--(A5) and $f_0=\infty$ and $f_\infty=0$, Theorem 5.7 has shown that there exists an unbounded
component $\mathscr{C}^+$ emanating from $(0,0)$ and joining to $(+\infty,+\infty)$.
Moreover, we also have the following theorem.
\\ \\
\textbf{Theorem 6.1.} \emph{Assume that} (A3)--(A5) \emph{hold and $f_0=\infty$ and $f_\infty=0$. Then, for any $\lambda\in(0,+\infty)$,
problem (\ref{pn}) has a unique positive solution $u_\lambda(x)$. Furthermore, such a solution $u_\lambda(x)$ satisfies the following properties:}
\\

%
(i) \emph{$u_\lambda(x)$ lies on $\mathscr{C}^+$;}

(ii) \emph{$u_\lambda(x)$ is continuous in $\lambda$, that is, if $\lambda\rightarrow \lambda_0$, then $\left\Vert u_\lambda-u_{\lambda_0}\right\Vert\rightarrow 0$.}
\\ \\
\textbf{Proof.} Suppose on the contrary that $(\lambda,u)$ and $(\lambda,v)$ are positive solutions satisfying $u\left(x_0\right)>v\left(x_0\right)$
at some point $x_0\in[0,T]$. We divide the rest proof into two cases.

{\it Case 1.} $u(0)\leq v(0)$.

In this case, it is not difficult to see that there is an interval $(c,d)$ such that $u>v$ in $(c,d)$ and $u(x)=v(x)$ at $x=c, d$.
By a direct computation one has
\begin{eqnarray}
-\int_c^d\left(\frac{u^p\varphi_p\left(v'\right)}{\varphi_p\left(v\right)}-u\varphi_p\left(u'\right)\right)'\,dx= \Gamma_1,\nonumber
\end{eqnarray}
where
\begin{eqnarray}
\Gamma_1=\int_c^d\left(\lambda m\left(\frac{f(v)}{\varphi_p(v)}-\frac{f(u)}{\varphi_p(u)}\right)u^p+\left(\left\vert u'\right\vert^p+(p-1)\left\vert \frac{uv'}{v}\right\vert^p-p\varphi_p\left(u\right)u'\varphi_p\left(\frac{v'}{v}\right)\right)\right)\,dx.\nonumber
\end{eqnarray}
It is not difficult to verify that
\begin{eqnarray}
\frac{u^p\varphi_p\left(v'\right)}{\varphi_p\left(v\right)}-u\varphi_p\left(u'\right)\geq 0\,\,\text{at}\,\, x=d,\nonumber
\end{eqnarray}
\begin{eqnarray}
\frac{u^p\varphi_p\left(v'\right)}{\varphi_p\left(v\right)}-u\varphi_p\left(u'\right)\leq0\,\,\text{at}\,\, x=c.\nonumber
\end{eqnarray}
Thus, $\Gamma_1\leq 0$. On the other hand, Young's inequality implies that
\begin{equation}
\left\vert u'\right\vert^p+(p-1)\left\vert \frac{uv'}{v}\right\vert^p-p\varphi_p\left(u\right)u'\varphi_p\left(\frac{v'}{v}\right)\geq 0.\nonumber
\end{equation}
In fact, we have
\begin{equation}
\left\vert u'\right\vert^p+(p-1)\left\vert \frac{uv'}{v}\right\vert^p-p\varphi_p\left(u\right)u'\varphi_p\left(\frac{v'}{v}\right)>0.\nonumber
\end{equation}
Indeed, if not, there exists $\mu\in \mathbb{R}$ such that $u=\mu v$ in $(c,d)$.
It is easy to show that $\mu=1$ since $u=v$ at $x=c, d$. This is a contradiction.
It follows that
\begin{eqnarray}
\int_c^d\lambda m\left(\frac{f(v)}{\varphi_p(v)}-\frac{f(u)}{\varphi_p(u)}\right)u^p\,dx<0,\nonumber
\end{eqnarray}
i.e.,
\begin{eqnarray}
\int_c^d\lambda m\left(\frac{f(v)}{\varphi_p(v)}\right)u^p\,dx< \int_c^d\lambda m\left(\frac{f(u)}{\varphi_p(u)}\right)u^p\,dx.\nonumber
\end{eqnarray}
Since $u(x)>v(x)$ in $(c,d)$, (A5) yields
\begin{eqnarray}
\int_c^d\lambda m f(u)u\,dx< \int_c^d\lambda m f(u)u\,dx.\nonumber
\end{eqnarray}
We get a contradiction.

\emph{Case 2.} $u(0)>v(0)$.

Obviously, there exists a positive $c_0\in(0,1)$ such that $v(0)=c_0 u(0)$. Let $\widetilde{u}=c_0u$. So we have $\widetilde{u}(0)=v(0)$, $\widetilde{u}(T)=v(T)$.
We consider the following three cases.

\emph{Case 2.1.} $c_0 u(x_0)>v(x_0)$.

Obviously, there is an interval $(e,f)$ such that $\widetilde{u}>v$ in $(e,f)$ and $\widetilde{u}(x)=v(x)$ at $x=e, f$.
By an argument similar to that of Case 1, we can obtain that
\begin{eqnarray}
\int_e^f\lambda m f(u)u\,dx< c_0^p\int_e^f\lambda m f(u)u\,dx<\int_e^f\lambda m f(u)u\,dx.\nonumber
\end{eqnarray}
This is a contradiction.

\emph{Case 2.2.} $c_0 u(x_0)<v(x_0)$.

Changing the roles of $\widetilde{u}$ and $v$ in the proof of Case 2.1, we can show that
\begin{eqnarray}
\int_g^h\lambda m f(v)v\,dx< \int_g^h\lambda m f(v)v\,dx,\nonumber
\end{eqnarray}
where $g,h\in[0,T]$ such that $\widetilde{u}>v$ in $(g,h)$ and $\widetilde{u}(x)=v(x)$ at $x=g, h$. We get a contradiction again.

\emph{Case 2.3.} $c_0 u(x_0)=v(x_0)$.

Firstly, we claim that $\widetilde{u}\not\equiv v$ on $[0,T]$. Suppose on the contrary that $\widetilde{u}\equiv v$ on $[0,T]$.
In view of equation (\ref{pn}), we can show that $c_0^{p-1}f(u)=f(v)$. By some simple computations, we obtain that
\begin{eqnarray}
\frac{f\left(\frac{v}{c_0}\right)}{\left(\frac{v}{c_0}\right)^{p-1}}=\frac{f(v)}{v^{p-1}}.\nonumber
\end{eqnarray}
Clearly, one has $v/c_0>v$ on $[0,T]$. It follows that
\begin{eqnarray}
\frac{f\left(\frac{v}{c_0}\right)}{\left(\frac{v}{c_0}\right)^{p-1}}<\frac{f(v)}{v^{p-1}}.\nonumber
\end{eqnarray}
We get a contradiction.

Thus, without loss of generality, we can assume that there exists $y_0\in\left(0,x_0\right)$ such that $\widetilde{u}\left(y_0\right)>v\left(y_0\right)$.
It is easy to see that there is an interval $(\alpha,\beta)$ such that $\widetilde{u}>v$ in $(\alpha,\beta)$ and $\widetilde{u}(x)=v(x)$ at $x=\alpha, \beta$.
Using an argument similar to that of Case 1, we have
\begin{eqnarray}
\int_\alpha^\beta\lambda m f(\widetilde{u})\widetilde{u}\,dx< \int_\alpha^\beta\lambda m f(\widetilde{u})\widetilde{u}\,dx,\nonumber
\end{eqnarray}
where $f(\widetilde{u})=c_0^{p-1}f(u)$. This is a contradiction.
Therefore, $\mathscr{C}^+$ is a curve.

%
Finally, we prove that $u_\lambda(x)$ is continuous with respect to $\lambda$. If $\lambda_0=0$, we define $u_0(x)\equiv 0$. It is obvious that
$\lim_{\lambda\rightarrow 0}\left\Vert u_\lambda\right\Vert=0$. Next, we assume that $\lambda_0>0$.
Let $\lambda>0$ such that $\lambda\rightarrow\lambda_0$ and $u_\lambda$ be the corresponding solutions.
Then we have that $u_\lambda$ is bounded since $\mathscr{C}^+$ does not blow up at a finite point.
By the compactness of $F_\lambda$, we obtain that for some convenient subsequence
$u_\lambda\rightarrow u$ in $X$. Clearly, we have $u=u_{\lambda_0}$. This completes the proof of the theorem.\qed\\

From Theorem 6.1, we can easily obtain the following result.
\\ \\
\textbf{Corollary 6.1.} \emph{Assume that} (A3)--(A5) \emph{hold. Then, for each $M\in(0,+\infty)$, there exists $\lambda_*\in(0,+\infty)$
such that (\ref{pn}) has a positive solution $u_*(x)$
with $\left\Vert u_*\right\Vert=M$.}
\\ \\
\textbf{Remark 6.1.} In [\ref{GKW}], the authors obtained the results similar to Theorem 6.1 under the assumptions of (A3), (A4) and
\\

(A6) $f:[0,+\infty)\rightarrow (0,+\infty)$ is nondecreasing, and there exists $\theta\in(0,1)$ such that
\begin{eqnarray}\label{6.1}
f(ks)\geq k^\theta f(s)\,\,\text{for}\,\, k\in(0,1)\,\,\text{and}\,\, u\in[0,+\infty).
\end{eqnarray}
Obviously, we do not need that $f$ to be nondecreasing in Theorem 6.1. In addition, (\ref{6.1}) implies the assumption of (A5) with $p=2$.
To see this, letting $0<s_1<s_2$, we show that $f\left(s_1\right)/s_1>f\left(s_2\right)/s_2$.
It is obvious that there exists a constant $k\in(0,1)$ such that $s_1=ks_2$.
Then we have
\begin{eqnarray}
\frac{f\left(s_2\right)}{s_2}\leq\frac{f\left(ks_2\right)}{k^\theta s_2}<\frac{f\left(ks_2\right)}{k s_2}\leq \frac{f\left(s_1\right)}{s_1}.\nonumber
\end{eqnarray}
Conversely, we cannot obtain (\ref{6.1}) from (A5) (with $p=2$), that is to say, (A5) is weaker than (\ref{6.1}) (with $p=2$). There is function
$f$ satisfying (A5) and not satisfying (\ref{6.1}). For example, let $f(s)=s^{\theta+\varepsilon}$ for some $\varepsilon\in(0,1-\theta)$. Clearly, $f$ satisfies (A5) (with $p=2$) but does not satisfy
(\ref{6.1}). Thus, our results have extended and improved the corresponding ones to [\ref{GKW}, Theorem 2.2] even in the case of $p=2$ in some sense.
\\ \\
\textbf{Remark 6.2.} In [\ref{GKW}], the authors also proved that
$u_\lambda$ is monotonic with respect to $\lambda$ in the case of $p=2$ and $q(x)=\rho^2$.
We conjecture that the solution $u_\lambda$ coming from Theorem 6.1 is also monotonic with respect to $\lambda$.

\section{One-sign solutions without signum condition}

\quad\, In Section 5, we have studied the existence of one-sign solutions for (\ref{pn}) under the signum condition.
Naturally, one may ask what will happen if $f$ does not satisfy signum condition.
In this section, we study problem (\ref{pn}) again but without signum condition.

\subsection{Unilateral global bifurcation from infinity}

\quad\, In this subsection, we study unilateral global bifurcation phenomena from infinity for problem (\ref{pb}).
Instead of (\ref{c1}), we assume that $g$ satisfies
\begin{equation}\label{c2}
\lim_{ \vert s\vert\rightarrow+\infty}\frac{g(x,s,\lambda)}{\vert s\vert^{p-1}}=0
\end{equation}
uniformly on $[0,T]$ and $\lambda$ on bounded sets.

We use $\mathscr{S}$ to denote the closure of the nontrivial solutions set of problem (\ref{pb}) in $\mathbb{R}\times X$.
We add the points $\{(\lambda,\infty)\big|\lambda\in \mathbb{R}\}$ to space $\mathbb{R}\times X$. Let $S_p$ denote the
spectral set of problem (\ref{pe}).\\

The main result of this subsection is the theorem below.
\\ \\
\textbf{Theorem 7.1.} \emph{Let the assumption (\ref{c2}) hold. There exists a component $\mathscr{D}_{\nu}^\sigma$ of
$\mathscr{\mathscr{S}}\cup\left(\lambda_0^\nu\times \{\infty\}\right)$, containing $\lambda_0^\nu\times \{\infty\}$.
Moreover if $\Lambda\subset \mathbb{R}$ is an interval such that $\Lambda\cap S_p=\lambda_0^\nu$ and $\mathscr{M}$ is
a neighborhood of $\lambda_0^\nu\times \{\infty\}$ whose projection on $\mathbb{R}$
lies in $\Lambda$ and whose projection on $X$ is bounded away from 0, then either}
\\

1$^o$. \emph{$\mathscr{D}_\nu^\sigma-\mathscr{M}$ is bounded in $\mathbb{R}\times X$ in which case $\mathscr{D}_\nu^\sigma-\mathscr{M}$
meets $\mathscr{R}=\{(\lambda,0)\big|\lambda\in \mathbb{R}\}$ or}
\\

2$^o$. \emph{$\mathscr{D}_\nu^\sigma-\mathscr{M}$ is unbounded.}\\

\emph{If 2$^o$ occurs and $\mathscr{D}_\nu^\sigma-\mathscr{M}$ has a bounded projection on $\mathbb{R}$, then
$\mathscr{D}_\nu^\sigma-\mathscr{M}$ meets $\widetilde{\lambda}\times \{\infty\}$ for some
$\widetilde{\lambda}\in S_p\setminus\left\{\lambda_0^+,\lambda_0^-\right\}$.}\\
\\
\textbf{Proof.} If $(\lambda,u)\in\mathscr{S}$ with $\Vert u\Vert\neq 0$, dividing (\ref{pb}) by $\Vert u\Vert^2$ and
setting $w=u/\Vert u\Vert^2$ yields
\begin{equation}\label{pnws}
\left\{
\begin{array}{l}
-\left(\varphi_p\left(w'\right)\right)'+q(x)\varphi_p(w)=\lambda m(x)\varphi_p(w)+\frac{g(x,u,\lambda)}{\Vert u\Vert^{2(p-1)}},\,\,0<x<T,\\
w(0)=w(T),\,\, w'(0)=w'(T).
\end{array}
\right.
\end{equation}
Define
\begin{equation}
\widetilde{g}(x,w,\lambda)=\left\{
\begin{array}{l}
\Vert w\Vert^{2(p-1)}g\left(x,\frac{w}{\Vert w\Vert^2},\lambda\right),\,\,\text{if}\,\, w\neq 0,\\
0,\,\,\,\,\,\quad\quad \quad\quad\quad\quad\quad\quad\quad\text{if}\,\, w=0.
\end{array}
\right.\nonumber
\end{equation}
Clearly, (\ref{pnws}) is equivalent to
\begin{equation}\label{pnws1}
\left\{
\begin{array}{l}
-\left(\varphi_p\left(w'\right)\right)'+q(x)\varphi_p(w)=\lambda m(x)\varphi_p(w)+\lambda m(x)\widetilde{g}(x,w,\lambda),\,\, x\in(0,T),\\
w(0)=w(T),\,\, w'(0)=w'(T).
\end{array}
\right.
\end{equation}
It is obvious that $(\lambda,0)$ is always the solution of (\ref{pnws1}).
By simple computation, we can show that the assumptions (\ref{c2}) implies
\begin{equation}
\widetilde{g}(x,s,\lambda)=o\left(\vert s\vert^{p-1}\right)\,\,\text{near}\,\, s=0,\,\,\text{uniformly for all}\,\,x\in[0,T]\,\,\text{and on bounded\,}\,\,
\lambda\,\,\text{intervals}.\nonumber
\end{equation}
Now applying Theorem 4.3 to problem (\ref{pnws1}), we have the component $\mathscr{C}_{\nu}^\sigma$ of
$\mathscr{S}\cup\left(\lambda_0^\nu\times \{0\}\right)$, containing $\lambda_0^\nu\times \{0\}$ is unbounded and lies in
$K^\sigma\cup\left(\lambda_0^\nu\times \{0\}\right)$.
Under the inversion $w\rightarrow w/\Vert w\Vert^2=u$, $\mathscr{C}_{\nu}^\sigma\rightarrow \mathscr{D}_\nu^\sigma$
satisfies problem (\ref{pb}). Clearly, $\mathscr{D}_\nu^\sigma$ satisfies the conclusions of this theorem.\qed\\

By Lemma 4.6 and the similar argument to prove [\ref{R0}, Corollary 1.8] with obvious changes,
we can obtain the following theorem.\\
\\
\textbf{Theorem 7.2.} \emph{There exists a neighborhood $\mathscr{N}\subset\mathscr{M}$ of $\lambda_0^\nu\times \{\infty\}$ such that
$(\lambda,u)\in\left(\mathscr{D}_\nu^\sigma\cap\mathscr{N}\right)\setminus\left\{\left(\lambda_0^\nu\times \{\infty\}\right)\right\}$ implies $(\lambda,u)
=\left(\lambda_0^\nu+o(1), \alpha \varphi_0^\nu+y\right)$
where $\Vert y\Vert=o(\vert \alpha\vert)$
at $\vert\alpha\vert=+\infty$.}

\subsection{Global behavior of the components of one-sign solutions}

\quad\, In this subsection, we study the problem (\ref{pn}) again but without signum condition. We only consider
the case of $f_0, f_\infty\in(0,+\infty)$. Other cases can be discussed similarly. The details are left to the reader.\\

Instead of (A2), we assume that $f$ satisfies:\\

(A7) there exist four constants $t_2\leq t_1<0<s_1\leq s_2$ such that $f\left(t_2\right)=f\left(t_1\right)
=f\left(s_2\right)=f\left(s_1\right)=f(0)=0$, $f(s)>0$ for
$s\in \left(t_2,t_1\right)\cup\left(0,s_1\right)\cup\left(s_2,+\infty\right)$ and $f(s)<0$ for
$s\in\left(-\infty,t_2\right)\cup\left(t_1,0\right)\cup\left(s_1,s_2\right)$.\\

We shall obtain the results similar to ones of [\ref{MXH}] for problem (\ref{pn}) in which the authors
only studied the existence of positive solutions with $p=2$. Note that
in the case of $p=2$, the authors of [\ref{MXH}] also
required that $f\in C^2(\mathbb{R},\mathbb{R})$ and satisfies $f''(s)<0$ for $s\in\left[0,s_1\right)$.
In this article, we drop these conditions completely. Hence, our results extend and improve the corresponding results of [\ref{MXH}].

Let $\xi$, $\eta\in C(\mathbb{R},\mathbb{R})$ be such that
\begin{equation}
f(u)=f_0\varphi_p(s)+\xi(s),\,\,\,f(s)=f_\infty\varphi_p(s)+\eta(s)\nonumber
\end{equation}
with
\begin{equation}
\lim_{\vert s\vert\rightarrow0}\frac{\xi(s)}{\varphi_p(s)}=0,\,\,\, \lim_{\vert s\vert\rightarrow+\infty}
\frac{\eta(s)}{\varphi_p(s)}=0.\nonumber
\end{equation}
Let us consider
\begin{equation}\label{pnws2}
\left\{
\begin{array}{l}
-\left(\varphi_p\left(u'\right)\right)'+q(x)\varphi_p(u)=\lambda f_0 m(x)\varphi_p(u)+\lambda m(x)\xi(u),\,\,0<x<T,\\
u(0)=u(T),\,\, u'(0)=u'(T)
\end{array}
\right.
\end{equation}
as a bifurcation problem from the trivial solution $u\equiv0$, and
\begin{equation}\label{pnws3}
\left\{
\begin{array}{l}
-\left(\varphi_p\left(u'\right)\right)'+q(x)\varphi_p(u)=\lambda f_\infty m(x)\varphi_p(u)+\lambda m(x)\eta(u),\,\,0<x<T,\\
u(0)=u(T),\,\, u'(0)=u'(T)
\end{array}
\right.
\end{equation}
as a bifurcation problem from infinity.

Applying Theorem 4.3 to (\ref{pnws2}), we have that there exists a
continuum $\mathscr{C}_\nu^\sigma$ of solutions to (\ref{pn}) joining $\left(\lambda_0^\nu/f_0,0\right)$ to infinity,
and $\left(\mathscr{C}_\nu^\sigma\setminus\left\{\left(\lambda_0^\nu/f_0,0\right)\right\}\right)\subseteq K^\sigma$.
Applying Theorem 7.1 to (\ref{pnws3}), we can show that there exists a
continuum $\mathscr{D}_\nu^\sigma$ of solutions to (\ref{pn}) meeting $\left(\lambda_0^\nu/f_\infty,\infty\right)$.
Moreover, Theorem 7.2 implies that $\left(\mathscr{D}_\nu^\sigma\setminus\left\{\left(\lambda_0^\nu/f_\infty,\infty\right)\right\}\right)\subseteq K^\sigma$.
\\

Next, we shall show that these two components are disjoint under the assumptions
(A7). Hence the essential role is played by the fact of whether $f$ possesses zeros in $\mathbb{R} \backslash\{0\}$.
\\ \\
\textbf{Theorem 7.3.} \emph{Let} (A7) \emph{hold. Then} \\

\indent (i) \emph{for $(\lambda,u)\in \left(\mathscr{C}_\nu^+\cup \mathscr{C}_\nu^-\right)$, we have that $t_1<u(x)<s_1$ for
all $x\in[0,T]$;}\\

(ii) \emph{for $(\lambda,u)\in \left(\mathscr{D}_\nu^+\cup \mathscr{D}_\nu^-\right)$, we have that either $\max_{x\in[0,T]}u(x)>s_2$
or $\min_{x\in[0,T]}u(x)<t_2$.}
\\ \\
\textbf{Proof.} We only prove for the case $(\lambda,u)\in\left(\mathscr{C}_+^+\cup \mathscr{D}_+^+\right)$
since the other cases can be proved similarly.
Suppose on the contrary that there exists
$(\lambda,u)\in\left(\mathscr{C}_+^+\cup \mathscr{D}_+^+\right)$ such that
either $\max\{u(x)\big|x\in[0,T]\}=s_1$
or $\min\{u(x)\big|x\in[0,T]\}=s_2$. We only treat the case of $\max\{u(x)\big|x\in[0,T]\}=s_1$
because the proof for the case of $\min\{u(x)\big|x\in[0,T]\}=s_2$ can be given similarly.

\emph{We claim that there exists $0<M<+\infty$ such that $f(s)\leq M\varphi_p\left(s_1-s\right)$ for any $s\in\left[0,s_1\right]$.}

Clearly, the claim is true for the case $s=0$ or $s=s_1$ by virtue of (A7). Suppose on the
contrary that there exists $s_0\in\left(0,s_1\right)$ such that
\begin{equation}
f\left(s_0\right)>M\varphi_p\left(s_1-s_0\right)\nonumber
\end{equation}
for any $M>0$. It follows that $M<f\left(s_0\right)/\varphi_p\left(s_1-s_0\right)$. This contradicts the arbitrary of $M$.

Now, let us consider the following problem
\begin{equation}
\left\{
\begin{array}{l}
-\left(\varphi_p\left(v'\right)\right)'+\lambda \vert m\vert M\varphi_p\left(v\right)
=\lambda \vert m\vert M\varphi_p\left(v\right)+q\varphi_p(u)-\lambda mf(u),\,\, x\in
(0,T),\\
s_1-u(0)\geq0,\,\,\,s_1-u(T)\geq0,
\end{array}
\right.\nonumber
\end{equation}
where $v=s_1-u$.
It is obvious that $f(s)\leq M\varphi_p\left(s_1-s\right)$ for any $s\in\left[0,s_1\right]$ implies
\begin{equation}
\left\{
\begin{array}{l}
-\left(\varphi_p\left(\left(s_1-u\right)'\right)\right)'+\lambda \vert m\vert M\varphi_p\left(s_1-u\right)\geq 0,\,\, x\in
(0,T),\\
s_1-u(0)\geq0,\,\,\,s_1-u(T)\geq0.
\end{array}
\right.\nonumber
\end{equation}
The strong maximum principle of [\ref{Mo}, Theorem 2] implies that $s_1 >u(x)$ on $(0,T)$.
Now, we show that $s_1 >u(x)$ on $[0,T]$. It suffices to show that $s_1-u(0)>0$ and $s_1-u(T)>0$.
Suppose on the contrary that $s_1-u(0)=0$ or $s_1-u(T)=0$, then we have $u(0)= u(T)=s_1$.
The strong maximum principle of [\ref{Mo}, Theorem 3] implies that $u'(0)<0$ and $u'(T)<0$ which
contradicts $u'(0)=u'(T)$. Thus, we have $s_1 >u(x)$ on $[0,T]$. This is a contradiction to
$\max\{u(x)\big|x\in[0,T]\}=s_1$.\qed
\\ \\
\textbf{Remark 7.1.} By Theorem 7.3, for any $(\lambda,u)\in \mathscr{C}_\nu^\sigma$, we can easily show
\begin{equation}
\Vert u\Vert_\infty<\max\left\{s_1,\left\vert t_1\right\vert\right\}:=s^*,\nonumber
\end{equation}
where $\Vert u\Vert_\infty=\max_{x\in[0,T]}\vert u\vert$.
Then we can easily deduce that
\begin{equation}
\Vert u\Vert\leq\left(\frac{\lambda T\Vert m\Vert_\infty\Vert u\Vert_\infty\max_{\vert s\vert
\leq \Vert u\Vert_\infty}\vert f(s)\vert}{p}\right)^{1/p}
\leq\left(\frac{\lambda T\Vert m\Vert_\infty s^*\max_{\vert s\vert\leq s^*}\vert f(s)\vert}{p}\right)^{1/p}.\nonumber
\end{equation}

In view of Theorem 7.1, Theorem 7.3 and Remark 7.1, using the similar argument to prove [\ref{M}, Corollary 2.1 and 2.2] with obvious changes,
we have the following two corollaries.
\\ \\
\textbf{Corollary 7.1.} \emph{Let} (A7) \emph{hold. Assume that $f_\infty>f_0$. Then}\\

(i) \emph{if $\lambda\in\left(\lambda_0^+/f_\infty,\lambda_0^+/f_0\right]\cup\left[\lambda_0^-/f_0,\lambda_0^-/f_\infty\right)$,
then (\ref{pn}) has at least two solutions
$u_{\infty}^+$ and
$u_{\infty}^-$, such that $u_{\infty}^+$ is positive and $u_{\infty}^-$ is negative;}

(ii) \emph{if $\lambda\in\left(\lambda_0^+/f_0,+\infty\right)\cup\left(-\infty,\lambda_0^-/f_0\right)$, then (\ref{pn})
has at least four solutions $u_{\infty}^+$,
$u_{\infty}^-$,
$u_{0}^+$ and $u_{0}^-$ such that $u_{\infty}^+$, $u_{0}^+$ are positive and $u_{\infty}^-$, $u_{0}^-$ are negative.}
\\ \\
\textbf{Corollary 7.2.} \emph{Let} (A7) \emph{hold. Assume that $f_0>f_\infty$. Then}\\

(i) \emph{if $\lambda\in\left(\lambda_0^+/f_0,\lambda_0^+/f_\infty\right]\cup\left[\lambda_0^-/f_\infty,\lambda_0^-/f_0\right)$,
then (\ref{pn}) has at least two solutions
$u_{0}^+$ and
$u_{0}^-$, such that $u_{0}^+$ is positive and $u_{0}^-$ is negative;}

(ii) \emph{if $\lambda\in\left(\lambda_0^+/f_\infty,+\infty\right)\cup\left(-\infty,\lambda_0^-/f_\infty\right)$, then (\ref{pn})
has at least four solutions $u_{\infty}^+$,
$u_{\infty}^-$,
$u_{0}^+$ and $u_{0}^-$ such that $u_{\infty}^+$, $u_{0}^+$ are positive and $u_{\infty}^-$, $u_{0}^-$ are negative.}
\\ \\
\textbf{Remark 7.2.} Besides the assumption of (A4) and (A7), we also assume that $f$ satisfies\\

(A8) $f(s)/\varphi_p(s)$ is strictly decreasing in $\left(0,s_1\right)$ and strictly increasing in $\left(t_1,0\right)$.\\

\noindent Then, by an argument similar to that of Theorem 6.1, we can show that
$\mathscr{C}_+^\sigma$ is a curve. Moreover, if $(\lambda,u_\lambda)\in \mathscr{C}_+^\nu$ then $u_\lambda$
is continuous in $\lambda$. We conjecture that this result is also valid if the assumption (A4) is removed.

\end{document}